\newif\ifpreprint
\preprinttrue 

\ifpreprint
	\documentclass[a4paper,oneside]{article}
	\usepackage{fullpage}
\fi

\usepackage{amsmath}
\usepackage{amssymb}

\ifpreprint
	\usepackage{amsthm}
	\usepackage{enumitem}
\fi
\usepackage{nicefrac}

\usepackage{graphicx}
\usepackage{subfig}
\usepackage{xcolor}
\definecolor{darkred}{RGB}{153,0,0}
\definecolor{darkblue}{RGB}{0,0,153}
\definecolor{darkgreen}{RGB}{0,153,0}

\ifpreprint
	\usepackage[colorlinks,
		citecolor=darkgreen,
		linkcolor=darkblue,
		urlcolor=darkred]{hyperref}
\fi

\ifpreprint
	\newcommand{\emphdef}[1]{\textsc{#1}}
\fi

\usepackage[ruled,vlined,linesnumbered,algosection,algo2e]{algorithm2e} 

\usepackage[capitalize,noabbrev,nameinlink]{cleveref} 

\ifpreprint
	\theoremstyle{plain}
\fi
\newtheorem{theorem}{Theorem}[section]
\newtheorem{lemma}[theorem]{Lemma}

\ifpreprint
	\theoremstyle{definition}
\fi
\newtheorem{assumption}[theorem]{Assumption}
\newtheorem{definition}[theorem]{Definition}
\ifpreprint
	\theoremstyle{remark}
\fi
\newtheorem{remark}[theorem]{Remark}
\newtheorem{examplex}[theorem]{Example}
\newenvironment{example}{\pushQED{\qed}\examplex}{\popQED\endexamplex}
\crefname{assumption}{Assumption}{Assumptions}

\ifpreprint
	\newcommand{\KeywordsAnd}{\and}
	\newcommand{\KeywordsEnd}{}
	\newcommand{\keywords}[1]{\par\noindent{\small\def\and{\unskip,\ }{\bf Keywords }#1.}\par}
	\newcommand{\subclass}[1]{\par\noindent{\small\def\and{\unskip,\ }{\bf AMS MSC }#1.}\par}
	\newcommand{\AMSand}{\and}
	\newcommand{\AMSend}{}
\fi

\crefname{line}{Step}{Steps} 
\SetKw{KwOr}{or}
\SetKw{KwElse}{else}
\SetKw{KwAnd}{and}
\SetKw{KwReturn}{return}
\SetKw{KwGoTo}{go to}

\newcommand{\emailLink}[1]{\textsc{email} \href{mailto:#1}{#1}}
\newcommand{\orcidLink}[1]{\textsc{orcid} \href{https://orcid.org/#1}{#1}}
\newcommand{\amsmscLink}[1]{\href{http://www.ams.org/mathscinet/msc/msc2020.html?t=#1}{#1}}

\newcommand{\N}{\mathbb{N}}
\newcommand{\Z}{\mathbb{Z}}
\newcommand{\R}{\mathbb{R}}

\DeclareMathOperator*{\argmin}{\arg\min}
\DeclareMathOperator{\dom}{dom}
\DeclareMathOperator{\proj}{proj}
\DeclareMathOperator{\prox}{prox}
\DeclareMathOperator{\dist}{dist}
\DeclareMathOperator{\indicator}{\delta}
\newcommand{\polarcone}{\circ}

\newcommand{\coloneqq}{:=}
\DeclareMathOperator*{\minimize}{minimize}
\DeclareMathOperator*{\maximize}{maximize}
\DeclareMathOperator{\stt}{subject~to}
\DeclareMathOperator{\wrt}{over}

\newcommand{\func}[3]{#1 \colon #2 \to #3}
\newcommand{\innprod}[2]{\langle #1, #2 \rangle}
\newcommand{\normalcone}{\mathcal{N}}
\newcommand{\limnormalcone}{\normalcone^{\textup{lim}}}
\newcommand{\conj}{\ast}
\newcommand{\jac}{\mathrm{J}}

\usepackage{tcolorbox}
\ifpreprint
	\newtcolorbox{mybox}[1][]{%
		left=0pt,
		right=0pt,
		top=0pt,
		bottom=0pt,
		colback=gray!12,
		colframe=gray!12,
		width=\dimexpr\textwidth\relax,
		enlarge left by=0mm,
		boxsep=5pt,
		arc=1pt,outer arc=1pt,
		#1
	}
	\newtcolorbox{mydefbox}[1][]{%
		left=0pt,
		right=0pt,
		top=0pt,
		bottom=0pt,
		colback=gray!12,
		colframe=gray!12,
		width=\dimexpr\textwidth\relax,
		enlarge left by=0mm,
		boxsep=5pt,
		arc=1pt,outer arc=1pt,
		#1
	}
\fi

\newcommand{\norm}[1]{\| #1 \|}
\newcommand{\closedball}{\mathbb{B}}
\newcommand{\lpfooter}{\textup{PL}}
\newcommand{\normlp}[1]{\norm{ #1 }_{\lpfooter}}
\newcommand{\lpball}{\closedball_{\lpfooter}}

\newcommand{\spaceX}{\R^n}
\newcommand{\spaceY}{\R^m}
\newcommand{\psimeas}{\Psi}
\newcommand{\lf}{\varphi}

\newcommand{\Lagr}{\mathcal{L}}
\newcommand{\feasmeas}{\mathcal{F}}
\newcommand{\XX}{\mathcal{X}}
\newcommand{\YY}{\mathcal{Y}}
\newcommand{\CC}{\mathcal{C}}
\newcommand{\KK}{\mathcal{K}}
\newcommand{\Ybounded}{\mathcal{Y}_s}
\newcommand{\II}{\mathcal{I}}
\newcommand{\dualLagr}{\mathcal{Q}}

\definecolor{mycolorNLP}{RGB}{27,158,119}%
\definecolor{mycolorCIA}{RGB}{217,95,2}%
\definecolor{mycolorMILA}{RGB}{117,112,179}%

\definecolor{color_turbo_al}{RGB}{228,26,28}%
\definecolor{color_turbo_dp}{RGB}{55,126,184}%
\definecolor{color_turbo_ip}{RGB}{77,175,74}%
\definecolor{color_turbo_guess}{RGB}{0,0,0}%

\newcommand{\TheAuthorADM}{Alberto De~Marchi}
\newcommand{\TheEmailADM}{alberto.demarchi@unibw.de}
\newcommand{\TheOrcidADM}{0000-0002-3545-6898}
\newcommand{\TheAffiliation}{%
	University of the Bundeswehr Munich,
	Department of Aerospace Engineering,
	Institute of Applied Mathematics and Scientific Computing,
	85577 Neubiberg, Germany%
}
\newcommand{\TheShortTitle}{Affordable mixed-integer Lagrangian methods}
\newcommand{\TheTitle}{\TheShortTitle: optimality conditions and convergence analysis}
\newcommand{\TheKeywords}{%
	Mixed-integer nonlinear programming\KeywordsAnd%
	Necessary optimality conditions\KeywordsAnd%
	Augmented Lagrangian framework\KeywordsAnd%
	Lagrangian duality\KeywordsAnd%
	Proximal point algorithm\KeywordsEnd%
}
\newcommand{\TheAMSsubj}{%
	\amsmscLink{65K05}\AMSand
	\amsmscLink{90C06}\AMSand
	\amsmscLink{90C11}\AMSand
	\amsmscLink{90C30}\AMSend
}
\newcommand{\TheAbstract}{%
	Necessary optimality conditions in Lagrangian form and the sequential minimization framework are
	extended to mixed-integer nonlinear optimization, without any convexity assumptions.
	Building upon a recently developed notion of local optimality for problems with
	polyhedral and integrality constraints, a characterization of local minimizers and critical points is
	given for problems including also nonlinear constraints.
	This approach lays the foundations for developing affordable sequential minimization algorithms with
	convergence guarantees to critical points from arbitrary initializations.
	A primal-dual perspective, a local saddle point property, and the dual relationships with the proximal point algorithm are also advanced in the presence of integer variables.
	Preliminary numerical results are presented for an augmented Lagrangian and an interior point method.
}%

\begin{document}

\ifpreprint
	\title{\bfseries \TheTitle}
	\author{\TheAuthorADM\thanks{\TheAffiliation.
			\emailLink{\TheEmailADM},
			\orcidLink{\TheOrcidADM}.}}
	\date{}

	\maketitle
	\begin{abstract}
		\TheAbstract
		
		\medskip
		
		\keywords{\TheKeywords}
		\subclass{\TheAMSsubj}
	\end{abstract}
	
	\tableofcontents
\fi


\section{Introduction}

Mixed-integer
nonlinear programming (MINLP) offers a versatile template for capturing a variety of tasks and applications,
but brings together ``the combinatorial difficulty of optimizing over discrete variable sets with the challenges of handling nonlinear functions'' \cite{belotti2013mixed}.
Originating from the integer programming community,
most approaches for MINLP rely on some sort of tree search
for seeking globally optimal solutions,
at least when some convexity is available.
Our focus is on affordable techniques
for addressing nonconvex MINLPs numerically.
Here, an optimization procedure is connotated as computationally ``affordable'' if it generates sequences globally convergent to points that satisfy some appropriate optimality conditions, though not necessarily to global minimizers.
In particular,
we are interested in iterative algorithms designed to converge to local solutions, in some sense, starting from arbitrary initial points \cite[Chapter 6]{birgin2014practical}.
Closely related to ``heuristics'' in the global optimization and integer programming community, these methods form the backbone of continuous optimization.
In practice, the potential benefit of reducing the explored search space is counteracted by weaker guarantees on the solution quality.
This tradeoff should allow us to handle large instances for a broad problem class,
but it requires defining a strong notion of local optimality,
with the aim of striking a balance between global but expensive minima and local but affordable critical points.

We seek a stationarity characterization that resembles, at least in spirit, the so called
Karush-Kuhn-Tucker (KKT) conditions in nonlinear programming (NLP).
Although ``in mixed-integer nonlinear programming, we do not know local optimality conditions comparable to the KKT conditions in continuous optimization'' \cite[Section 2]{exler2012comparative},
some advancements have been made based on an excess of multipliers and separation theorems \cite{jahn2018lagrange}.
In an attempt to upgrade our understanding, we study here a criticality concept for nonconvex MINLPs in simple Lagrangian terms.
Building upon the \emph{partial localization} approach and the corresponding optimality notions developed in \cite{demarchi2025mixed} for simply constrained problems,
we dedicate this work to characterizing ``local'' minima with a Lagrangian perspective and then establishing convergence results for a class of augmented Lagrangian (AL) methods.

A mixed-integer linearization algorithm (MILA) was proposed in \cite{demarchi2025mixed} to address the minimization of a smooth function over a feasible set with mixed-integer linear structure, namely MINLP without nonlinear constraints.
Even beyond AL schemes,
we are motivated by the sequential (partially) unconstrained minimization framework \cite{fiacco1968nonlinear},
which includes (shifted) penalty \cite{birgin2014practical} and barrier (or interior point) methods \cite{waechter2006implementation},
to handle nonlinear constraints while taking advantage of the affordable solver of \cite{demarchi2025mixed} for tackling the subproblems.
The present work provides solid theoretical foundations for this algorithmic design paradigm, exemplified by AL methods.
We discuss how this framework can be used to design other algorithms for MINLP,
and in particular we indicate how similar arguments apply also to interior point approaches on the line of \cite{demarchi2024interior}.
Methods based on sequential mixed-integer quadratic programming \cite{exler2012comparative,quirynen2021sequential} could benefit from these theoretical advances too.
Other numerical approaches for MINLP, such as global methods or decomposition techniques \cite{belotti2013mixed,sager2011combinatorial},
could also exploit these principled heuristics to refine initial guesses, generate tighter bounds, and promote faster convergence.

Beyond numerical methods for MINLP,
we enrich the theoretical framework and first-order analysis of mixed-integer optimization in Lagrangian terms,
inspired by the celebrated KKT conditions in nonlinear programming.
In the spirit of \cite{jahn2018lagrange,steck2018lagrange,rockafellar2023convergence},
we develop a theory of KKT-critical points, complemented by Lagrangian duality, saddle point properties, and relationships with the proximal points algorithm.

The problem template with nonconvex smooth objective and polyhedral, integrality, and nonlinear set-membership constraints reads
\begin{align}
	\tag{P}\label{eq:P}
	\minimize~ f(x) \quad
	\wrt~ x\in\XX \quad
	\stt~ c(x) \in \CC
	,
\end{align}
where $x\in\XX\subset\spaceX$ are decision variables,
$\func{f}{\XX}{\R}$ and $\func{c}{\XX}{\spaceY}$ are continuously differentiable functions,
$\CC\subset\spaceY$ is a nonempty closed convex set (projection-friendly in practice),
and $\XX$ is a nonempty closed set with mixed-integer linear structure \cite{demarchi2025mixed}.
In particular, set $\XX$ admits a description in the form of intersection between a closed convex polyhedral set $\overline{\XX}\subseteq\spaceX$ (that is, finitely many linear inequalities) and integrality constraints defined by some index set $\II\subset\{1,2,\ldots,n\}$:
\begin{equation*}
	\XX
	\coloneqq
	\overline{\XX}
	\cap
	\left\{
	x \in \spaceX
	\,\middle|\,
	x_i \in \Z ~~ \forall i\in\II
	\right\}
	.
\end{equation*}
In the following,
we may refer to a partition of decision variables $x$ into real-valued and integer-valued ones,
respectively $\{x_i \,|\, i\notin\II\}$ and $\{x_i \,|\, i\in\II\}$.
Furthermore, patterning \cite{demarchi2025mixed}, we consider the following blanket assumptions.

\begin{mybox}
	\begin{assumption}\label[assumption]{ass:P}
		With regard to \eqref{eq:P},
		\begin{enumerate}[label=(\alph*),series=assumptions,ref=\alph*]
			\item\label{ass:P:wellposed}
			$\inf \left\{ f(x) \,|\, x\in\XX ,\, c(x)\in\CC \right\} \in \R$;
			\item\label{ass:P:locLipschitzDiff}%
			functions $f$ and $c$ are continuously differentiable;
			\item\label{ass:P:integerBounded}%
			for all $i\in\II$ the set $\left\{ a \in\Z \,|\, x\in\XX ,\, x_i=a \right\}$ is bounded.
		\end{enumerate}
	\end{assumption}
\end{mybox}
The basic \cref{ass:P}\eqref{ass:P:wellposed} ensures that \eqref{eq:P} is well-posed, namely that it is feasible and a solution exists; it is adopted in the theoretical analysis and it is \emph{not} needed for the proposed algorithm to operate. Practical solvers typically include algorithmic safeguards and mechanisms to detect infeasibility or unboundedness and return with appropriate warnings.
Differentiability of $f$ and $c$ in \cref{ass:P}\eqref{ass:P:locLipschitzDiff} is intended with respect to real- and integer-valued variables,
treating them all as real-valued ones to avoid exotic definitions or approximations,
such as those in \cite{exler2012comparative}.
A practical situation that satisfies \cref{ass:P}\eqref{ass:P:locLipschitzDiff} is when $f$ and $c$ depend linearly on the integer-valued variables,
as supposed in \cite{quirynen2021sequential}.
Finally, \cref{ass:P}\eqref{ass:P:integerBounded} guarantees that admissible values (with respect to $\XX$ alone) for the integer-valued decision variables lie in a bounded set.
As it applies to integer-valued variables only,
this boundedness requirement is reasonable and often satisfied in practice (trivially for binary variables).
Following \cite{demarchi2025mixed},
we take advantage of \cref{ass:P}\eqref{ass:P:integerBounded} to construct compact neighborhoods without explicitly localizing the integer-valued components.

\subsection{Prompt, Outline and Contribution}

A major motivation for this work is the application to optimal control of hybrid dynamical systems, whose (time discretized) models comprise real- and integer-valued variables, nonlinear possibly nonsmooth dynamics, and combinatorial constraints.
Of particular interest is the case of mixed-integer optimal control, where the time structure has been exploited to design decomposition methods with approximation guarantees \cite{sager2011combinatorial}.
Relying on relaxation and subsequent combinatorial integral approximation (CIA),
this strategy exploits mature technology for NLP and mixed-integer linear programming (MILP), as well as the peculiar structure of optimal control problems \cite{buerger2020pycombina}.
However, since the classical CIA does not take into account the system dynamics nor path constraints, it can generate infeasible trajectories.
Moreover, when combinatorial constraints are present (such as dwell time constraints), the CIA sub-optimality bounds might be severely affected \cite{zeile2021mixed}.
To overcome these issues, recent works \cite{buerger2023gauss,ghezzi2023voronoi} have proposed to formulate the CIA problem as a mixed-integer quadratic program (MIQP) that locally approximates the MINLP of interest.

In the same spirit, we advocate here for preserving the structure of \eqref{eq:P} as much as possible, while seeking good quality, not necessarily global, solutions.
This avenue was explored numerically in \cite{nikitina2025hybrid} and it is further motivated here with an example presented in \cref{sec:example_mila_cia}, where a
direct comparison on a simple problem illustrates the advantages of holding on to the integrality constraint, without relaxing it.
Animated by the numerical approach proposed in \cite{demarchi2025mixed} and the extensions foreseen there,
we build the theoretical foundations for sequential minimization algorithms to address \eqref{eq:P},
establishing convergence results under suitable assumptions.
Our monolithic strategy provides convergence guarantees and can be adopted as a framework to combine several techniques, such as relaxations, integral approximations and feasibility pumps \cite{dambrosio2012storm,belotti2013mixed}.

Our contributions can be summarized as follows:
\begin{itemize}
	\item We derive and analyse necessary optimality conditions for \eqref{eq:P} in Lagrangian form, comparable to the KKT system in continuous optimization---see \cref{sec:StationarityConcepts}.
	\item We prove the global convergence of a safeguarded augmented Lagrangian algorithm, providing a solid theoretical support for generalizing the affordable approach of \cite{demarchi2025mixed} to sequential minimization schemes for MINLP---see \cref{sec:AugLagFramework}.
	\item The Lagrangian system is further characterized in primal-dual terms, recovering saddle-point properties and a dual relationship with the proximal point algorithm---see \cref{sec:FurtherCharacterizations}.
\end{itemize}
The main goal of this work is to lay solid theoretical foundations that
support the numerical approach of \cite{nikitina2025hybrid} and
provide the basis for further methodological developments.
Although comprehensive computational investigations are beyond the scope of this paper, some numerical results showcased in \cref{sec:NumericalResults} substantiate the proposed algorithmic framework.

\subsection{Notation and Preliminaries}
The set of natural, integer, and real numbers are denoted by $\N$, $\Z$, $\R$.
The appearing spaces are equipped with the standard Euclidean inner product $\innprod{\cdot}{\cdot}$ and norm $\| \cdot \|$.
Given a nonempty subset $\CC$ of $\spaceY$,
the \emphdef{indicator} $\func{\indicator_{\CC}}{\spaceY}{\R\cup\{\infty\}}$,
the \emphdef{projection} $\func{\proj_{\CC}}{\spaceY}{\CC}$,
and the \emphdef{distance} $\func{\dist_{\CC}}{\spaceY}{\R}$
are defined respectively by
\begin{align*}
	\indicator_{\CC}(v)
	\coloneqq{}&
	\begin{cases}
		0 & \text{if}~v\in\CC, \\
		\infty & \text{otherwise},
	\end{cases} &
	\proj_{\CC}(v)
	\coloneqq{}&
	\argmin_{z\in\CC} \| z - v \| , &
	\dist_{\CC}(v)
	\coloneqq{}&
	\min_{z\in\CC} \| z - v \| .
\end{align*}
The \emphdef{normal cone} $\normalcone_{\CC}(z)$ of set $\CC \subseteq\spaceY$ at $z\in\CC$ is given by
\begin{equation*}
	\normalcone_{\CC}(z) \coloneqq \left\{ v\in\spaceY \,\middle\vert\, \forall u\in\CC\colon \innprod{v}{u-z} \leq 0 \right\} .
\end{equation*}
For formal completeness, we define $\normalcone_{\CC}(z) \coloneqq \emptyset$ if $z\notin \CC$.
We will make use of the following well-known characterizations
valid for a closed convex set $\CC \subseteq \spaceY$ \cite{bauschke2017convex}:
\begin{equation}
	\label{eq:projCharacterization}
	u\in\proj_{\CC}(z)
	\iff
	\forall w\in \CC \colon~ \innprod{z-u}{w-u} \leq 0
	,
\end{equation}
\begin{equation}
	\label{eq:normalConeCharacterization}
	u\in\normalcone_{\CC}(z)
	\iff
	\forall \alpha > 0 \colon~ z = \proj_{\CC}(z + \alpha u)
	\iff
	\exists \alpha > 0 \colon~ z = \proj_{\CC}(z + \alpha u)
	.
\end{equation}

\section{Optimality Concepts}
\label{sec:OptimalityConcepts}

A point $\bar{x} \in \spaceX$ is called \emphdef{feasible} for \eqref{eq:P} if it satisfies the constraints there, namely $\bar{x}\in\XX$ and $c(\bar{x}) \in \CC$.
It is also clear how to define a \emph{global} solution, or minimizer, $x^\star$ for \eqref{eq:P}: a feasible point where the optimal objective value is attained, namely
\begin{align*}
	x^\star \in{} \XX ,\quad
	c(x^\star) \in{} \CC ,\quad
	\forall x\in\XX, c(x)\in{} \CC \colon~ f(x^\star) \leq f(x) .
\end{align*}
But then, what constitutes a suitable notion of \emph{local} minimizer?

Answers to this important question affect not only the quality of what we refer to as ``solutions'',
but they do influence also the design of numerical methods.
Before handling nonlinear constraints with the Lagrangian formalism, let us review the approach proposed in \cite{demarchi2025mixed} for simply constrained problems.

\subsection{Neighborhoods with Partial Localization}

Local notions, as opposed to global ones, depend on the concept of neighborhood and this, in turn, is very delicate in the mixed-integer context of \eqref{eq:P}.
Following \cite{demarchi2025mixed},
we denote by $\normlp{\cdot}$ an operator mapping $x$ into a norm of the real-valued entries of $x$, that is, given the index set $\II$.
Prominent examples are $\normlp{v} \coloneqq \max_{i\notin\II} |v_i|$
and $\normlp{v} \coloneqq \sum_{i\notin\II} |v_i|$,
associated with $\ell_\infty$ and $\ell_1$ norms, respectively.
The notation ``PL'' stands for \emphdef{partial localization}, owing to the fact that PL-balls
\begin{equation*}
	\lpball(x,\Delta)
	\coloneqq
	\left\{
	w\in\spaceX \,|\, \normlp{w-x}\leq\Delta
	\right\}
\end{equation*}
identify a neighborhood for the real-valued components and not for the integer-valued ones, which remain free.
For this reason, PL-balls are \emph{not} compact sets in general.
Nevertheless,
the intersection
$\XX \cap \lpball(x,\Delta)$
is always a compact set,
thanks to \cref{ass:P}\eqref{ass:P:integerBounded},
and thus represents a reasonable neighborhood of $x$
---and a valid trust region stipulation---
for any $x\in\XX$ and $\Delta\geq 0$.
Before proceeding,
we should mention that adopting a polyhedral norm to define $\normlp{\cdot}$ is favourable in practice,
as the mixed-integer \emph{linear} structure is not lost in the subproblems,
but the theory applies with any norm.

A local concept of solution for \eqref{eq:P} can now be defined by means of these (partial) neighborhoods.
Inspired by \cite[Definition~2]{demarchi2025mixed}, local and global minimizers for \eqref{eq:P} are characterized as follows.
\begin{mydefbox}
	\begin{definition}\label[definition]{def:minimizer}
		A point $\bar{x}\in\spaceX$ is called a \emphdef{local minimizer} for \eqref{eq:P} if it is feasible
		and
		there exists $\Delta>0$ such that
		$f(\bar{x}) \leq f(x)$ for all feasible $x\in\lpball(\bar{x},\Delta)$.
		If the latter property additionally holds for all $\Delta>0$, then $\bar{x}$ is called a \emphdef{global minimizer}.
	\end{definition}
\end{mydefbox}
For instances of \eqref{eq:P} without integer-valued variables, namely $\II\coloneqq\emptyset$,
\cref{def:minimizer} recovers the classical notion of local minima in nonlinear programming.
Conversely, without real-valued variables, namely $\II\coloneqq\{1,2,\ldots,n\}$,
\eqref{eq:P} is an integer program and \cref{def:minimizer} effectively requires a global solution (since there is no actual localization in this case).
Thus, we can observe that monitoring neighborhoods with $\normlp{\cdot}$ leads to a stronger local optimality concept than a plain adaptation of continuous notions into the mixed-integer realm---for instance, using Euclidean neighborhoods in $\R^n$.
In fact, local minimizers in the sense of \cref{def:minimizer} are also stronger than those obtained by
`fixing the integer variables and optimizing over the continuous ones',
since a certificate of local optimality must consider all feasible points in $\lpball(\bar{x},\Delta)$,
which may contain several integer configurations.
Conversely, the combinatorial structure in \eqref{eq:P} should be simple enough for practical purposes, e.g., mixed-integer linear.

Before delving into KKT-like optimality conditions for \eqref{eq:P},
let us recall some solution concepts for problems without nonlinear constraints.
Following \cite{demarchi2025mixed}, consider the minimization of $\func{\lf}{\XX}{\R}$ over $\XX$ as a basic template:
\begin{equation}
	\label{eq:Punc}
	\minimize~ \lf(x) \qquad
	\wrt~ x\in\XX .
\end{equation}
A local notion of solutions for \eqref{eq:Punc} is proposed in \cite[Definition~2]{demarchi2025mixed}, inspired by \cite[Definition~3.1]{byrd2005convergence} for the analogous minimization over a \emph{convex} set.
A first-order optimality measure associated to \eqref{eq:Punc} (that is, to function $\lf$ and set $\XX$) is defined in \cite[Equation~4]{demarchi2025mixed} and provides a metric $\psimeas_{\lf,\XX}$ to monitor ``optimality'':
for all $x\in\XX$ and $\Delta>0$ it is given by
\begin{equation}
	\label{eq:psimeas}
	\psimeas_{\lf,\XX}(x, \Delta)
	\coloneqq
	\max_{w\in\XX\cap\lpball(x,\Delta)} \innprod{\nabla \lf(x)}{x - w}
	\geq
	0
	.
\end{equation}
Since $x, w\in\XX$ in \eqref{eq:psimeas}, $\psimeas_{\lf,\XX}(\cdot,\Delta)$ is bounded from below by zero for all $\Delta>0$.
Note that the maximization in \eqref{eq:psimeas} is over all feasible points in $\lpball(x,\Delta)$.
Then, a first-order optimality concept for the ``simply constrained'' problem \eqref{eq:Punc} is defined as follows; cf. \cite[Definition~3]{demarchi2025mixed}.
\begin{mydefbox}
	\begin{definition}\label[definition]{def:criticality}
		Given some $\varepsilon>0$ and $\Delta > 0$,
		a point $\bar{x} \in \spaceX$ is called \emphdef{$\varepsilon$-$\Delta$-critical} for \eqref{eq:Punc} if $\bar{x} \in \XX$ and
		$\psimeas_{\lf,\XX}(\bar{x}, \Delta) \leq \varepsilon$.
		Given some $\varepsilon>0$,
		a point $\bar{x} \in \spaceX$ is called \emphdef{$\varepsilon$-critical} for \eqref{eq:Punc} if it is $\varepsilon$-$\Delta$-critical for some $\Delta>0$.
		A $0$-critical point is simply called \emphdef{critical}.
	\end{definition}
\end{mydefbox}
\cref{def:criticality} provides a valid concept to characterize candidate minimizers, necessary for optimality,
which is stronger than plain (M-)\emph{stationarity}; see
\cite[Section~2.2]{demarchi2025mixed} and \cref{example:optimality_notions} below.
The criticality notion for ``unconstrained'', or simply constrained, problems \eqref{eq:Punc} will become important to characterize solutions to intermediate, auxiliary problems (referred to as subproblems).
Moreover, defining an approximate counterpart of criticality allows us to consider inexact subproblem solutions, a strategy often (if not always) adopted in sequential minimization methods.
This is useful in accommodating iterative subsolvers with asymptotic convergence, and then in exploiting this property to reduce the overall computational effort.

\medskip

Since the quality of subproblem solutions eventually affects the (outer loop) iterates,
stronger optimality notions can lead to better performance,
as illustrated with the following example.
It turns out that, in the context of mixed-integer problems,
criticality based on PL-balls provides not only good candidates for minimizers,
but often it is also easier to compute than projection-based continuous counterparts.

\begin{example}[Optimality, criticality and stationarity]\label{example:optimality_notions}
	Consider a two-dimensional problem of the form \eqref{eq:Punc} with decision variable $x\coloneqq (u,z)$:
	\begin{equation}\label{eq:toy_example_criticality}
		\minimize_{u,\, z}\quad u^2 \qquad\stt\quad z \leq u \leq 1+z ,\quad z\in\{0,1\}.
	\end{equation}
	The feasible set $\XX$ is the union of two line segments, as depicted in \cref{fig:toy_example_criticality:base},
	and the global minimizer for \eqref{eq:toy_example_criticality} is $x^\star\coloneqq (0,0)$, with objective $f^\star=0$.
	The characterization of our focus point $\bar{x} \coloneqq (1,1)$ is open to debate.
	With $f(\bar{x})=1>f^\star$, it is clearly not a global solution,
	but whether it is a ``local'' minimizer or not depends on the point-of-view.
	
	\begin{itemize}
	\item
	Continuous optimization:
	\begin{itemize}
		\item From a variational analysis perspective, $\bar{x}$ is a feasible stationary point for \eqref{eq:toy_example_criticality},
		since the inclusion $0 \in \nabla f(\bar{x})+\limnormalcone_\XX(\bar{x})$ is valid.
		Moreover, in view of
		\cite[Definition 3.1, Proposition 3.5]{themelis2018forward},
		the point $\bar{x}$ is not only stationary but also critical for \eqref{eq:toy_example_criticality},
		since $\bar{x} \in \proj_\XX(\bar{x} - \gamma \nabla f(\bar{x}))$ holds for all stepsizes $\gamma \in (0,\nicefrac{1}{4}]$.
		\item
		Even considering a connected enlargement of the feasible set,
		as in \cref{fig:toy_example_criticality:enlarged},
		$\bar{x}$ is locally optimal in the classical sense of continuous optimization (that is, taking a ball around $\bar{x}$ in $\R^2$).
		In fact, most (if not all) nonlinear programming solvers initialized at $\bar{x}$ would stop there, declaring a successful solve,
		since $\bar{x}$ is B-stationary (that is, there are no feasible first-order directions of descent).
	\end{itemize}
	\item
	Heuristics for MINLP:
	\begin{itemize}
		\item Fixing the integer variable at $z\coloneqq 1$, the value $f(\bar{x})=1$ at $\bar{u}\coloneqq 1$ is optimal.
		Moreover, fixing the continuous variable at $u\coloneqq 1$, the value $f(\bar{x})=1$ at $\bar{z}\coloneqq 1$ is also optimal,
		so that $\bar{x}\coloneqq (1,1)$ is reasonably deemed ``locally'' optimal.
	\end{itemize}
	\item
	Neighborhoods with partial localization:
	\begin{itemize}
		\item Given any $\Delta>0$, the corresponding PL-neighborhood of $\bar{x}$ is given by two disconnected line segments
		and covers a point $x_\Delta$ with better objective value than that of $\bar{x}$.
		For instance, with $\Delta\in(0,1]$, we find $f(x_\Delta) = (1-\Delta)^2 < f(\bar{x})$ at $x_\Delta\coloneqq(1-\Delta,0) \in \XX\cap\lpball(\bar{x},\Delta)$.
		Therefore, according to \cref{def:criticality}, the point $\bar{x}$ is \emph{not} critical for \eqref{eq:toy_example_criticality} and, in particular,
		$\bar{x}$ is \emph{not} a local minimizer in the sense of \cref{def:minimizer}.
	\end{itemize}
\end{itemize}
	
	\begin{figure}
		\centering%
		\subfloat[][\emph{Two-dimensional toy problem}\label{fig:toy_example_criticality:base}]
		{%
		\includegraphics{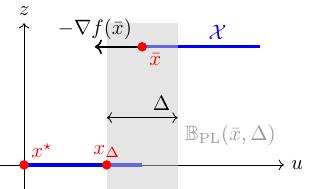}
		}
		\hfill
		\subfloat[][\emph{$\ldots$ with enlarged feasible set}\label{fig:toy_example_criticality:enlarged}]
		{%
		\includegraphics{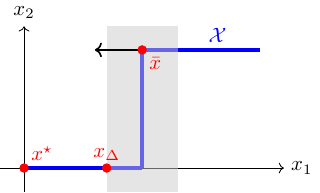}
		}
		\hfill
		\subfloat[][\emph{$\ldots$ with reduced feasible set}\label{fig:toy_example_criticality:gap}]
		{%
		\includegraphics{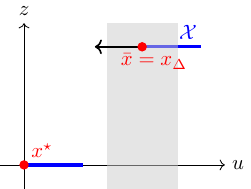}
		}
		\caption{%
			Illustration of two-dimensional toy problems discussed in \cref{example:optimality_notions}.
			Each panel depicts a feasible set $\XX$ (thick blue line), the global solution $x^\star \coloneqq (0,0)$,
			the focus point $\bar{x} \coloneqq (1,1)$ and the (negative) gradient there $\nabla f(\bar{x})$,
			the PL-ball $\lpball(\bar{x},\Delta)$ with radius $\Delta>0$ (shaded gray area),
			and the trust-region update $x_\Delta\in\XX\cap\lpball(\bar{x},\Delta)$.
			While $\bar{x}$ is a stationary point in all three cases,
			it can be deemed critical (in the PL sense) only in scenario \ref{fig:toy_example_criticality:gap} (for some sufficiently small $\Delta>0$).
			In scenarios \ref{fig:toy_example_criticality:base} and \ref{fig:toy_example_criticality:enlarged},
			sufficient decrease is attained by taking the step from $\bar{x}$ to $x_\Delta$.
			Thus, the PL-based criticality of \cref{def:criticality} is stronger than stationarity.
		}%
		\label{fig:toy_example_criticality}
	\end{figure}
	
	These observations pertain the characterization of ``solutions'' but have also computational consequences:
	methods of projected-gradient type may detect ``local optimality'' at $\bar{x}$ and stop there,
	whereas methods (such as MILA) based on PL-neighborhoods deem $\bar{x}$ unsatisfactory and continue their search process
	(discovering the segment with $z=0$ and improving the objective value until they find the global minimizer $x^\star$).
	Effectively, PL-based methods easily escape the spurious point $\bar{x}$ where the others above can get trapped
	since they are oblivious of the discrete structure and rely on a continuous view.

	An example where the point $\bar{x}$ is also critical (in the PL sense) is illustrated in \cref{fig:toy_example_criticality:gap},
	with the feasible set $\XX$ representing the constraints $z \leq u \leq \nicefrac{1}{2}+z$ and $z\in\{0,1\}$.
	Particularly, the criticality measure is $\psimeas(\bar{x},\Delta)=0$ for all $\Delta\in(0,\nicefrac{1}{2})$, and it is $x_\Delta=\bar{x}$.

	Although there are no guarantees that PL-based methods will always outperform the others in terms of objective value,
	we can expect them to deliver ``better'' solutions as they are based on a (strictly) stronger optimality notion.
	In practice, this means that they \emph{could} make further progress where other methods stop, as illustrated by the toy problem \eqref{eq:toy_example_criticality} above, while the reverse situation is not possible.
\end{example}

Before extending the concept of critical points to the more general problem \eqref{eq:P}, a clarification on the role of $\Delta$ is in order.
\begin{remark}\label{rem:radius_only_theory}
By \cref{def:criticality}, an $\varepsilon$-critical point $\bar{x}$ is associated to \emph{some} radius $\Delta>0$.
The same happens when characterizing critical points in nonsmooth nonconvex optimization, where a (proximal) stepsize effectively acts as the radius $\Delta$ here; cf. \cite[Definition 3.1(ii)]{themelis2018forward}.
However, since the value of $\Delta$ need not be known, the mixed-integer Lagrangian framework developed here for \eqref{eq:P} is \emph{not} restricted, or specific, to the MILA of \cite{demarchi2025mixed}.
This fact is witnessed by \cref{alg:AbstractALM,alg:AbstractIP} below, where only an approximate critical point is required from the subsolver and there is no mention of $\Delta$.
In principle, projected-gradient and Frank-Wolfe methods could also be adopted as subsolvers.
In practice, though, the availability of affordable subsolvers appears limited: Frank-Wolfe schemes often rely on some convexity in the problem \cite{hendrych2022convex}, whereas Euclidean projections lead to MIQPs in general, in contrast to MILPs arising from \eqref{eq:psimeas}, hindering the performance of projected schemes.
\end{remark}

\subsection{Stationarity Concepts and Lagrangian Analysis}
\label{sec:StationarityConcepts}

What is a ``critical point'' for \eqref{eq:P}?
Treating the nonlinear constraints explicitly,
let the Lagrangian function $\func{\Lagr}{\XX\times\spaceY}{\R}$ associated to \eqref{eq:P} be defined, as usual, by
\begin{equation}
	\label{eq:classicalLagrangian}
	\Lagr(x,y)
	\coloneqq
	f(x) + \innprod{y}{c(x)}
	.
\end{equation}
From the viewpoint of nonlinear programming,
where stationarity of the Lagrangian plays a crucial role,
we consider the following notion for KKT-like points of \eqref{eq:P} based on \cref{def:criticality}.
Then, we are going to establish the (asymptotic) necessity of KKT-criticality for local optimality.
Related concepts and results can be found in \cite{demarchi2024interior,demarchi2023constrained,demarchi2024local,rockafellar2023convergence}.
\begin{mydefbox}
	\begin{definition}\label[definition]{def:KKTcritical}
		Given some $\Delta > 0$,
		a point $\bar{x} \in \spaceX$ is called \emphdef{$\Delta$-KKT-critical} for \eqref{eq:P} if $\bar{x}\in\XX$ and there exists a multiplier $y \in \spaceY$ such that
		\begin{equation*}
			\psimeas_{\Lagr(\cdot,y),\XX}(\bar{x},\Delta) = 0 
			\quad\text{and}\quad
			y \in \normalcone_\CC(c(\bar{x})) .
		\end{equation*}
		A point $\bar{x} \in \spaceX$ is called \emphdef{KKT-critical} for \eqref{eq:P} if it is $\Delta$-KKT-critical for some $\Delta > 0$.
	\end{definition}
\end{mydefbox}
KKT-criticality implicitly requires feasibility, since the normal cone $\normalcone_\CC(c(\bar{x}))$ must be nonempty.
Moreover, by \eqref{eq:psimeas} the first condition can be rewritten as
\begin{equation*}
	\min_{x\in\XX\cap\lpball(\bar{x},\Delta)} \innprod{\nabla f(\bar{x}) + \jac c(\bar{x})^\top y}{x - \bar{x}}
	=
	0 ,
\end{equation*}
meaning that the Lagrangian function cannot be (locally) further minimized with respect to $x$ while maintaining mixed-integer linear feasibility, in the sense of \cref{def:criticality},
effectively replacing stationarity with criticality.%
\footnote{%
	Although unclear whether multipliers can be affine sensitivities or not in MINLP, \cref{def:KKTcritical}'s introduction of multipliers $y$ for \eqref{eq:P} is harmless because they are associated to classical constraints only, which are smooth by \cref{ass:P}\eqref{ass:P:locLipschitzDiff}.
	This observation is supported by the role played by multipliers $y$ in the proof of \cref{thm:KKTnecessaryCondition}.
}
An asymptotic counterpart of \cref{def:KKTcritical} (also referred to as sequential or approximate)
proves to be a key tool for convergence analysis; cf. \cite[Definition~3.1]{birgin2014practical},
\cite{demarchi2023constrained}.
\begin{mydefbox}
	\begin{definition}\label[definition]{def:AKKTcritical}
		A point $\bar{x} \in \spaceX$ is called \emphdef{asymptotically KKT-critical} (AKKT) for \eqref{eq:P} if $\bar{x} \in \XX$ and there exist sequences $\{x^k\}\subset\spaceX$, $\{y^k\}\subset\spaceY$, $\{z^k\}\subseteq\CC$, and $\{\Delta_k\}\subset \R_{++}$ such that $x^k \to \bar{x}$ and
		\begin{equation*}
			\psimeas_{\Lagr(\cdot,y^k),\XX}(x^k,\Delta_k) \to 0 ,\quad
			y^k \in \normalcone_\CC( z^k ) ,\quad
			c(x^k) - z^k \to 0 .
		\end{equation*}
	\end{definition}
\end{mydefbox}
If a sequence $\{x^k\}$ has an accumulation point which is AKKT-critical,
then finite termination can be attained with an approximate KKT-critical point, for any given tolerance $\varepsilon>0$.
\begin{mydefbox}
	\begin{definition}\label[definition]{def:epsKKTcritical}
		Given some $\varepsilon\geq0$,
		a point $\bar{x} \in \spaceX$ is called \emphdef{$\varepsilon$-KKT-critical} for \eqref{eq:P} if $\bar{x} \in \XX$ and there exist a multiplier $y\in\spaceY$, a vector $z\in\CC$, and some $\Delta>0$ such that
		\begin{equation*}
			\psimeas_{\Lagr(\cdot,y),\XX}(\bar{x},\Delta) \leq \varepsilon ,\quad
			y \in \normalcone_\CC( z ) ,\quad
			\| c(\bar{x}) - z \| \leq \varepsilon .
		\end{equation*}
		A $0$-KKT-critical point is simply called \emphdef{KKT-critical}.
	\end{definition}
\end{mydefbox}

We can now establish a link between minimizers in the sense of \cref{def:minimizer} and KKT-like critical points.
A local minimizer for \eqref{eq:P} is
KKT-critical under validity of a suitable qualification condition.
However, each local minimizer of \eqref{eq:P} is always AKKT-critical, regardless of additional regularity.
Related results can be found in \cite{birgin2014practical,demarchi2023constrained}.
\begin{mybox}
	\begin{theorem}\label[theorem]{thm:KKTnecessaryCondition}
		Let $x^\star\in\spaceX$ be a local minimizer for \eqref{eq:P}.
		Then, $x^\star$ is AKKT-critical.
	\end{theorem}
\end{mybox}
\begin{proof}\label{proof:KKTnecessaryCondition}
	By local optimality of $x^\star$ for \eqref{eq:P} there exists $\delta > 0$ such that $f(x^\star) \leq f(x)$ is valid for all feasible $x \in \lpball(x^\star,\delta)$; cf. \cref{def:minimizer}.
	Consequently, $x^\star$ is the unique global minimizer of the localized problem
	\begin{align}
		\minimize~&f(x) + \|x-x^\star\|^2 \qquad
		\wrt~x \in\XX \cap \lpball(x^\star,\delta)
		\label{eq:proof:localizedProblem2}\\
		\stt~&c(x) \in \CC . \nonumber
	\end{align}
	Slightly deviating from the proof of \cite[Proposition~2.5]{demarchi2023constrained},
	let us consider the penalized surrogate problem
	\begin{equation}
		\minimize~\pi_k(x) \qquad \wrt~x \in\XX \cap \lpball(x^\star,\delta) ,
		\label{eq:proof:penaltySurrogateProblem2}
	\end{equation}
	where
	\begin{equation*}
		\pi_k(x) \coloneqq f(x) + \|x-x^\star\|^2 + \rho_k \dist_{\CC}^2( c(x) ),
	\end{equation*}
	$k \in \N$ is arbitrary, $\rho_k > 0$, and the sequence $\{\rho_k\}_{k\in\N}$ satisfies $\rho_k\to\infty$ as $k \to\infty$.
	
	Noting that the objective function of this optimization problem is lower semicontinuous while its feasible set is nonempty and compact (by feasibility of $x^\star$, trust region stipulation, and \cref{ass:P}\eqref{ass:P:integerBounded}),
	it possesses a global minimizer $x^k \in \XX$ for each $k \in \N$,
	owing to Weierstrass' extreme value theorem.
	Without loss of generality, we assume $x^k \to \widetilde{x}$ for some $\widetilde{x}\in\XX \cap \lpball(x^\star,\delta)$.
	
	We now argue that $\widetilde{x} = x^\star$.
	To this end, we note that $x^\star$ is feasible to \eqref{eq:proof:penaltySurrogateProblem2} with $c(x^\star) \in \CC$,
	which yields for each $k\in\N$ the (uniform, upper) estimate
	\begin{equation}
		\label{eq:proof:penaltyUpperEstimate2}
		\pi_k(x^k)
		=
		f(x^k) + \|x^k - x^\star\|^2 + \rho_k \dist_{\CC}^2(c(x^k)) \leq f(x^\star)
		.
	\end{equation}
	Using $\rho_k \to \infty$, lower semicontinuity of $f$, finiteness of $f(x^\star)$, closedness of $\CC$, and the convergence $c(x^k) \to c(\widetilde{x})$, taking the limit for $k\to\infty$ in \eqref{eq:proof:penaltyUpperEstimate2} gives $c(\widetilde{x}) \in \CC$.
	Therefore, $\widetilde{x}$ is feasible for \eqref{eq:P} and local optimality of $x^\star$ for \eqref{eq:P} implies $f(x^\star) \leq f(\widetilde{x})$.
	Furthermore, exploiting \eqref{eq:proof:penaltyUpperEstimate2} and the optimality of each $x^k \in \XX$, we find
	\begin{equation*}
		f(\widetilde{x}) + \|\widetilde{x} - x^\star\|^2
		\leq{}
		\liminf_{k\to\infty} \pi_k(x^k)
		\leq{}
		f(x^\star)
		\leq
		f(\widetilde{x})
		.
	\end{equation*}
	Hence, $\widetilde{x} = x^\star$.
	Now we may assume without loss of generality that $\{x^k\}$ is taken from the interior of $\lpball(x^\star,\delta)$, as this is eventually the case, since $x^k\to x^\star$.
	Thus, for each $k \in\N$, $x^k$ globally minimizes $\pi_k$ over $\XX$,
	see \eqref{eq:proof:penaltySurrogateProblem2},
	whose relevant criticality condition (necessary for optimality \cite[Proposition~1]{demarchi2025mixed}) reads, for some $\Delta_k>0$,
	\begin{align*}
		0
		={}&
		\psimeas_{\pi_k,\XX}(x^k,\Delta_k)
		={}
		\max_{w\in\XX\cap\lpball(x^k,\Delta_k)} \innprod{\nabla_x \Lagr(x^k,y^k) + 2 (x^k - x^\star)}{x^k - w}
	\end{align*}
	where we set $y^k \coloneqq 2 \rho_k [c(x^k) - \proj_{\CC}(c(x^k)) ]$ for each $k \in\N$.
	Now, owing to continuous differentiability of $\Lagr$ and compactness of $\XX\cap\lpball(x^k,\Delta_k)$,
	by $x^k\to x^\star\in\XX$ we have
	\begin{align*}
		\lim_{k\to\infty}
		\psimeas_{\Lagr(\cdot,y^k),\XX}(x^k,\Delta_k)
		={}&
		\lim_{k\to\infty}
		\max_{w\in\XX\cap\lpball(x^k,\Delta_k)} \innprod{\nabla_x \Lagr(x^k,y^k)}{x^k - w }
		={}
		\lim_{k\to\infty}
		\psimeas_{\pi_k,\XX}(x^k,\Delta_k)
		={}
		0 .
	\end{align*}
	Thus, the conditions in \cref{def:AKKTcritical} are a consequence of $x^k\to x^\star$. Overall, this shows that any local minimizer $x^\star$ for \eqref{eq:P} is AKKT-critical.
\end{proof}

Bridging the gap between AKKT- and KKT-criticality
requires some sort of constraint qualifications (CQ), such as the well-known LICQ and MFCQ.
In general, these are geometric conditions or stability properties that bound the set of Lagrange multipliers and thus guarantee that local minimizers are indeed KKT-critical;
see \cite{birgin2014practical} for a more detailed discussion.

\section{Augmented Lagrangian Framework}
\label{sec:AugLagFramework}

Let us consider \eqref{eq:P} under \cref{ass:P},
which, under the lens of continuous optimization,
can be seen as a nonlinear program with mixed-integer linear constraints.
Since the restriction to $\XX$ is nonrelaxable but easy to satisfy, in the sense that we treat it as hard while assuming that the associated MILPs are efficiently solved,
such constraint can be treated in a way essentially different from how nonlinear constraints are handled \cite{andreani2008augmented,birgin2014practical}.

The algorithms examined in the following are sequential minimization schemes designed to generate iterates whose limit points are AKKT-critical for \eqref{eq:P}
and so candidate minimizers, according to \cref{thm:KKTnecessaryCondition}.
\cref{alg:AbstractALM,alg:AbstractIP} below are implemented relying on the (approximate) PL-based criticality concept of \cref{def:criticality},
even though they could stand on mere stationarity.
We make this algorithmic choice explicit to highlight how the theoretical notion of local optimality illustrated with \cref{example:optimality_notions} can benefit numerical practice,
since the quality of limit points depends on how well each subproblem is solved.
Loosely writing, the stronger the criticality notion adopted, the greater the chances of finding high-quality minimizers for \eqref{eq:P}.

In the following \cref{sec:Algorithm} we study an AL method as an epitome for the class of sequential minimization schemes \cite{fiacco1968nonlinear}.
A theoretical characterization of the abstract \cref{alg:AbstractALM} is detailed in \cref{sec:ConvergenceAnalysis}, and the adjustments needed when considering other sequential minimization schemes (such a barrier methods) are sketched in \cref{sec:OtherSchemes}.

\subsection{Algorithm}
\label{sec:Algorithm}

The AL framework has been broadly investigated and developed, giving rise to a variety of multifaceted ideas, of which we only scratch the surface here.
The interested reader may refer to
\cite{birgin2014practical} for an overview,
to \cite{demarchi2024local,rockafellar2023convergence,steck2018lagrange} for theoretical advances,
and to \cite{demarchi2023constrained,sopasakis2020open} for numerical aspects.
The main ingredient of AL methods is the AL function
$\func{\Lagr_\mu}{\XX\times\spaceY}{\R}$,
whose definition associated to \eqref{eq:P} is
\begin{equation}
	\label{eq:augLagr}
	\Lagr_\mu(x,y)
	\coloneqq{}
	f(x) + \frac{1}{2 \mu} \dist_{\CC}^2\left( c(x) + \mu y \right) - \frac{\mu}{2} \|y\|^2
\end{equation}
for some penalty parameter $\mu > 0$ and multiplier estimate $y\in\spaceY$.
This is a \emph{partial} AL function in that it does not relax the simple constraint $x\in\XX$, which is kept explicit in each subproblem.
Notice that $\Lagr$ and $\Lagr_\mu$ are smooth, with respect to both, primal and dual variables $x$ and $y$, thanks to \cref{ass:P}\eqref{ass:P:locLipschitzDiff} and convexity of $\CC$.
For later use,
the partial derivatives of $\Lagr_\mu$ read
\begin{align}\label{eq:ALderivatives}
	\nabla_x \Lagr_\mu(x,y)
	={}&
	\nabla f(x) + \jac c(x)^\top y_\mu(x,y) ,&
	\nabla_y \Lagr_\mu(x,y)
	={}&
	c(x) - s_\mu(x,y)
\end{align}
where
\begin{align}\label{eq:ALauxiliaries}
	s_\mu(x,y) \coloneqq{}& \proj_{\CC}(c(x) + \mu y) ,&
	y_\mu(x,y) \coloneqq{}& y + \frac{c(x) - s_\mu(x,y)}{\mu} .
\end{align}
Following the basic pattern of AL methods,
\cref{alg:AbstractALM} proceeds by minimizing the AL function at each iteration, possibly inexactly and up to criticality, and updating the multiplier estimates and penalty parameters \cite[Section~4.1]{birgin2014practical}.
Augmented Lagrangian subproblems require to
\begin{equation}
	\label{eq:ALsubproblem}
	\minimize~ \Lagr_\mu(x,\widehat{y}) \quad
	\wrt~ x \in\XX
\end{equation}
given some $\mu>0$ and $\widehat{y}\in\spaceY$.%
\footnote{%
	It should be stressed that, within the scope of this paper, subproblem \eqref{eq:ALsubproblem} is indeed easier than the original \eqref{eq:P}.
	Since it has only mixed-integer linear constraints, it can be tackled with the \emph{local} approach of \cite{demarchi2025mixed}.
	To be sure, seeking a local solution to \eqref{eq:P}, there is no need to employ global techniques (such as spatial branch-and-bound, among others) to find a global solution for \eqref{eq:ALsubproblem}, making it relatively practical to solve \eqref{eq:ALsubproblem} up to (approximate) criticality.
}
Feasibility of \eqref{eq:ALsubproblem} follows from $\XX$ being nonempty, whereas well-posedness is due to (lower semi)continuity of $\Lagr_\mu(\cdot,\widehat{y})$ and is guaranteed if, e.g., $\XX$ is compact or $f$ is bounded from below in $\XX$.
In fact, the existence of subproblem solutions is often just assumed, see, e.g., \cite[Assumption 6.1]{birgin2014practical}.
Algorithmically, this difficulty could be circumvented by
complementing the AL subproblems \eqref{eq:ALsubproblem} with a localizing constraint, e.g., of trust region type
\cite[Remark~5.1]{demarchi2024implicit}.
However, 
as for the original problem \eqref{eq:P}, whose solutions exist according to \cref{ass:P}\eqref{ass:P:wellposed},
we merely assume that all subproblems are well-posed.
Analogous in spirit to prox-boundedness \cite[Definition~1.23]{rockafellar2009variational},
our \cref{ass:P:proxbounded}
is weaker than typical coercivity or (level) boundedness assumptions
but sufficient to yield well-posed subproblems.

\begin{mybox}
	\begin{assumption}\label[assumption]{ass:P:proxbounded}
		With regard to \eqref{eq:P} and \cref{alg:AbstractALM},
		there exists $\bar{\mu}>0$ such that
		for all $\mu\in(0,\bar{\mu}]$ and $\widehat{y}\in\Ybounded$
		the function $\Lagr_\mu(\cdot,\widehat{y})$ is bounded from below over $\XX$.
	\end{assumption}
\end{mybox}

This allows us to focus on the mixed-integer extension of generic AL methods to address MINLP.
A practical implementation of the solver should provide mechanisms for detecting infeasibility and unboundedness, as discussed in \cite{dambrosio2012storm,quirynen2021sequential}.

\begin{algorithm2e}[tbh]
	\DontPrintSemicolon
	Select $\mu_0\in(0,\bar{\mu}]$, $\varepsilon_0, \eta_0 > 0$,
	$\kappa_\mu, \theta_\mu \in (0,1)$, and $\Ybounded \subseteq \spaceY$ bounded\;
	\For{$j = 0,1,2\ldots$}{
		Select $\widehat{y}^j \in \Ybounded$\label{step:AbstractALM:ysafe}\;
		Find an $\varepsilon_j$-critical point $x^j$ for $\Lagr_{\mu_j}(\cdot,\widehat{y}^j)$ over $\XX$\label{step:AbstractALM:subproblem}\tcp*{subproblem}
		Set $z^j \gets \proj_{\CC}(c(x^j) + \mu_j \widehat{y}^j)$, $v^j \gets c(x^j) - z^j$, and $y^j \gets \widehat{y}^j + \mu_j^{-1} v^j$ \label{step:AbstractALM:y}\label{step:AbstractALM:cviol}\;
		\If{$j=0$ \KwOr $\|v^j\| \leq \max\{\eta_j,\theta_\mu \|v^{j-1}\|\}$\label{step:AbstractALM:check_mu}}{%
			set $\mu_{j+1} \gets \mu_j$, \KwElse select $\mu_{j+1} \in (0,\kappa_\mu \mu_j]$\label{step:AbstractALM:update_mu}
		}
		Select $\varepsilon_{j+1}, \eta_{j+1} \geq 0$ such that $\{\varepsilon_j\}, \{\eta_j\}\to 0$\label{step:AbstractALM:innertol}\;
	}
	\caption{Abstract safeguarded augmented Lagrangian method for \eqref{eq:P}}
	\label{alg:AbstractALM}
\end{algorithm2e}

The scheme outlined in \cref{alg:AbstractALM} is often referred to as \emph{safeguarded}
because the multiplier estimates $\widehat{y}$ are not allowed to grow too fast compared to the penalty parameter $\mu$ \cite{birgin2014practical,steck2018lagrange,sopasakis2020open,demarchi2023constrained}.
In particular, it is required that $\|\mu_j \widehat{y}^j\|\to 0$ as $\mu_j\to0$,
so that stronger global convergence properties can be attained.
As a simple mechanism to ensure this property,
multiplier estimates $\widehat{y}$ in \cref{alg:AbstractALM} are drawn from a bounded set $\Ybounded\subseteq\spaceY$.
The dual safeguarding set $\Ybounded$ can be a generic hyperbox or can be tailored to the constraint set $\CC$ at hand \cite{sopasakis2020open}---see \cref{sec:LagrangianDuality} below.

Subproblems \eqref{eq:ALsubproblem} can be solved up to approximate criticality:
given $\varepsilon_j$, at \cref{step:AbstractALM:subproblem} we seek an $\varepsilon_j$-critical point $x^j \in \XX$ for $\Lagr_{\mu_j}(\cdot,\widehat{y}^j)$,
in the sense of \cref{def:criticality}.
For this task one can employ the mixed-integer linearization algorithm of \cite{demarchi2025mixed}, with guarantee of finite termination under \cref{ass:P,ass:P:proxbounded}.
Although the trust region radius $\Delta_j$ associated to the $\varepsilon_j$-criticality certificate does not need to be computed,
it will be considered formally for the theoretical analysis; cf. \cref{rem:radius_only_theory}.
Given a (possibly inexact, first-order) solution $x$ to \eqref{eq:ALsubproblem}, the dual update rule at \cref{step:AbstractALM:y}
is designed toward the identity
\begin{equation}
	\label{eq:gradLagrIdentity}
	\nabla_x \Lagr_\mu(x,\widehat{y})
	={}
	\nabla f(x) + \jac c(x)^\top y
	={}
	\nabla_x \Lagr(x,y) ,
\end{equation}
as usual in AL methods.
This allows to monitor the (outer) convergence with the (inner) subproblem tolerance; cf. \cref{lem:iterates} below.

Finally, \cref{step:AbstractALM:check_mu,step:AbstractALM:update_mu,step:AbstractALM:innertol} are dedicated to monitoring primal feasibility (namely the conditions involving $z^k$ in \cref{def:AKKTcritical}) and updating the penalty parameter $\mu$ accordingly.
Note that considering a sequence of primal tolerances $\{\eta_j\}$ allows to monitor primal convergence from a global perspective,
slightly relaxing in fact other classical update rules \cite{birgin2014practical,demarchi2024implicit}.

\subsection{Convergence Analysis}
\label{sec:ConvergenceAnalysis}

\Cref{alg:AbstractALM} belongs to the family of safeguarded AL schemes \cite{steck2018lagrange} and, by keeping the mixed-integer linear constraints explicit in subproblem \eqref{eq:ALsubproblem}, as opposed to relaxing them,
it closely resembles the AL scheme with lower-level constraints of \cite{andreani2008augmented,birgin2014practical}.
Thus, the following proofs pattern those found in classical AL literature,
but they all have the peculiarity of dealing with some trust region radius $\Delta$.
This feature is due to the deliberate choice of (approximate) criticality over mere stationarity when solving \eqref{eq:ALsubproblem} at \cref{step:AbstractALM:subproblem}, leading to stronger optimality notions and, plausibly, better solutions; see the discussion in \cref{example:optimality_notions}.

We begin our asymptotic analysis by collecting useful properties to characterize the iterations generated by \cref{alg:AbstractALM}.

\begin{mybox}
	\begin{lemma}\label[lemma]{lem:iterates}
		Let \cref{ass:P,ass:P:proxbounded} hold for \eqref{eq:P} and consider the iterates of \cref{alg:AbstractALM}.
		Then, for each $j\in\N$,
		\cref{step:AbstractALM:subproblem} is well-posed and the iterates satisfy
		$x^j\in\XX$, $z^j\in\CC$, $y^j \in \normalcone_{\CC}(z^j)$, $\nabla_x \Lagr_{\mu_j}(x^j,\widehat{y}^j) = \nabla_x \Lagr(x^j,y^j)$, and there exists some $\Delta_j>0$ such that $\psimeas_{\Lagr(\cdot,y^j),\XX}(x^j,\Delta_j) \leq \varepsilon_j$.
	\end{lemma}
\end{mybox}
\begin{proof}
	Well-definedness of \cref{alg:AbstractALM} follows from the existence of solutions to the AL subproblems,
	which in turn is due to the standing \cref{ass:P,ass:P:proxbounded}.
	In particular, the feasible set $\XX$ is nonempty and closed, and the continuous real-valued cost function $\Lagr_{\mu_j}(\cdot,\widehat{y}^j)$ is lower bound over $\XX$, since $\mu_j \leq \bar{\mu}$, for all $j\in\N$.
	
	Then, it is apparent that $x^j\in\XX$ and $z^j\in\CC$ for each $j\in\N$.
	Moreover, the assignments at \cref{step:AbstractALM:y} gives that
	$z^j
	\coloneqq
	\proj_{\CC}(c(x^j) + \mu_j \widehat{y}^j)
	=
	c(x^j) + \mu_j \widehat{y}^j - \mu_j y^j$,
	which is equivalent to $y^j\in\normalcone_{\CC}(z^j)$ by \eqref{eq:normalConeCharacterization} and convexity of $\CC$.
	By construction \eqref{eq:gradLagrIdentity}, the dual update rule
	readily yields $\nabla_x \Lagr_{\mu_j}(x^j,\widehat{y}^j) = \nabla_x \Lagr(x^j,y^j)$,
	and so the upper bound on the criticality measure and the existence of a suitable $\Delta_j$ follow from \cref{step:AbstractALM:subproblem}.
\end{proof}

We now turn to investigating properties of accumulation points, assuming their existence (which may follow from coercivity or level boundedness arguments).
The following convergence results for \cref{alg:AbstractALM} provides fundamental theoretical support for the numerical approach envisioned in \cite{demarchi2025mixed} to deal with nonlinear constraints, based on \cite{fiacco1968nonlinear}.
With \cref{thm:globalConvergence} we establish that feasible accumulation points of $\{x^j\}$ are AKKT-critical;
see \cite[Thm 3.3]{demarchi2023constrained}, \cite[Thm 3.6]{demarchi2024implicit} for analogous results.

\begin{mybox}
	\begin{theorem}\label[theorem]{thm:globalConvergence}
		Let \cref{ass:P,ass:P:proxbounded} hold.
		Consider a sequence $\{x^j\}$ generated by \cref{alg:AbstractALM}.
		Let $x^\star$ be an accumulation point of $\{x^j\}$ and $\{x^j\}_{j\in J}$ a subsequence such that $x^j \to_J x^\star$.
		If $x^\star$ is feasible for \eqref{eq:P}, then $x^\star$ is AKKT-critical for \eqref{eq:P}.
	\end{theorem}
\end{mybox}
\begin{proof}
	It is implicitly assumed \cref{alg:AbstractALM} generates an infinite sequence of iterates $\{x^j\}$ with accumulation point $x^\star$.
	Now we claim that the subsequences $\{x^j\}_{j\in J}$, $\{y^j\}_{j\in J}$, $\{z^j\}_{j\in J}$, $\{\Delta^j\}_{j\in J}$ satisfy the properties in \cref{def:AKKTcritical},
	thus showing that $x^\star$ is AKKT-critical for \eqref{eq:P}.
	From \eqref{eq:psimeas} and \cref{lem:iterates} we have that for all $j\in\N$
	\begin{equation*}
		0
		\leq
		\psimeas_{\Lagr(\cdot,y^j),\XX}(x^j,\Delta_j)
		\leq
		\varepsilon_j
	\end{equation*}
	for some $\Delta_j>0$.
	Hence, dual feasibility holds asymptotically owing to $\varepsilon_j \to 0$.
	
	By assumption we have $x^j \to_J x^\star$ with $x^\star$ feasible for \eqref{eq:P}, namely $x^\star\in\XX$ and $c(x^\star) \in \CC$.
	\cref{lem:iterates} implies also that $y^j \in \normalcone_\CC(z^j)$ for each $j\in\N$.
	Finally, to demonstrate that $c(x^j) - z^j \to_J 0$ we consider two cases:
	\begin{itemize}
		\item If $\{\mu_j\}$ is bounded away from zero, the conditions at \cref{step:AbstractALM:update_mu,step:AbstractALM:check_mu} of \cref{alg:AbstractALM} and the construction of $\{\eta_j\}$ imply that $\|v^j\| \coloneqq \| c(x^j) - z^j \| \to 0$, hence the assertion.
		\item If $\mu_j \to 0$, we exploit continuity of $c$, boundedness of $\{\widehat{y}^j\} \subseteq \Ybounded$, feasibility of $x^\star$, and closedness of $\CC$.
		Combining these properties gives $c(x^j) + \mu_j \widehat{y}^j \to_J c(x^\star) \in \CC$ as $x^j \to_J x^\star$.
		Therefore, $z^j \to_J c(x^\star)$ as well, hence $c(x^j)-z^j\to_J 0$.
	\end{itemize}
	Overall, this proves that $x^\star$ is AKKT-critical for \eqref{eq:P}.
\end{proof}

In contrast with global methods \cite[Chapter~5]{birgin2014practical}, \cite[Section~4.2]{demarchi2024local},
adopting affordable solvers for addressing \eqref{eq:ALsubproblem} at \cref{step:AbstractALM:subproblem}
impedes to guarantee that, in general, accumulation points are feasible or (globally) minimize an infeasibility measure.
Thus, despite feasibility granted by \cref{ass:P}\eqref{ass:P:wellposed}, \cref{alg:AbstractALM} may not approach feasible points.
In practice, however, 
for any fixed $\mu>0$ and $\widehat{y}\in\spaceY$,
the AL subproblem \eqref{eq:ALsubproblem} is equivalent to
\begin{equation*}
	\label{eq:ALsubproblemScaled}
	\minimize~
	\mu f(x) + \frac{1}{2} \dist_{\CC}^2\left( c(x) + \mu \widehat{y} \right) \quad
	\wrt~ x \in\XX
	.
\end{equation*}
Hence, one can expect to find at least critical points of an infeasibility measure,
as attested by the following result.
Notice that this property requires mere boundedness of $\{\varepsilon_j\}$; cf. \cite[Thm~6.3]{birgin2014practical},
\cite[Proposition~3.7]{demarchi2024implicit}.

\begin{mybox}
	\begin{theorem}\label[theorem]{thm:minimalInfeasibility}
		Let \cref{ass:P,ass:P:proxbounded} hold.
		Consider a sequence $\{x^j\}$ generated by \cref{alg:AbstractALM} with $\{\varepsilon_j\}$ merely bounded.
		Let $x^\star$ be an accumulation point of $\{x^j\}$ and $\{x^j\}_{j\in J}$ a subsequence such that $x^j \to_J x^\star$.
		Then, $x^\star$ is a critical point for the feasibility problem
		\begin{equation*}
			\minimize~
			\feasmeas(x) \coloneqq \frac{1}{2}\dist_{\CC}^2 ( c(x) ) \quad
			\wrt~ x\in\XX .
		\end{equation*}
	\end{theorem}
\end{mybox}
\begin{proof}\label{proof:minimalInfeasibility}
	It is implicitly assumed that \cref{alg:AbstractALM} generates an infinite sequence of iterates $\{x^j\}$ with accumulation point $x^\star$.
	If $\{\mu_j\}$ is bounded away from zero,
	the conditions at \cref{step:AbstractALM:update_mu,step:AbstractALM:check_mu} of \cref{alg:AbstractALM} and the construction of $\{\eta_j\}$ imply that $\|v^j\| \coloneqq \| c(x^j) - z^j \| \to 0$.
	By the upper bound $\|v^j\| \geq \dist_{\CC}( c(x^j) )$ for each $j\in\N$, since $z^j\in\CC$,
	taking the limit $j\to\infty$ yields $c(x^\star)\in\CC$ by continuity.
	Then, since $x^j\in\XX$ for all $j\in\N$ and $\XX$ is closed, $x^\star$ is feasible for \eqref{eq:P}.
	Thus, $x^\star$ is a global minimizer of the feasibility problem and, by continuous differentiability of the objective function therein, $x^\star$ is critical for the feasibility problem.
	
	Let us focus now on the case where $\{\mu_j\} \searrow 0$ and $x^\star\in\XX$ is infeasible for \eqref{eq:P}.
	First, we express what criticality entails for the feasibility problem above:
	a point $\bar{x}\in\spaceX$ is critical if $\bar{x}\in\XX$ and there exists some $\Delta > 0$ such that $\psimeas_{\feasmeas,\XX}(\bar{x}, \Delta) = 0$.
	Now, owing to \eqref{eq:psimeas} and \cref{step:AbstractALM:subproblem}, for all $j\in\N$ it is
	\begin{equation*}
		\varepsilon_j
		\geq
		\psimeas_{\Lagr_{\mu_j}(\cdot,\widehat{y}^j),\XX}(x^j, \Delta_j)
		=
		\max_{w\in\XX\cap\lpball(x^j,\Delta_j)}
		\innprod{\nabla_x \Lagr_{\mu_j}(x^j,\widehat{y}^j)}{x^j - w}
		\geq
		0 .
	\end{equation*}
	Multiplying by $\mu_j>0$, by boundedness of $\{\varepsilon_j\}$ we have
	\begin{equation*}
		0
		\leq
		\max_{w\in\XX\cap\lpball(x^j,\Delta_j)}
		\innprod{\mu_j \nabla_x \Lagr_{\mu_j}(x^j,\widehat{y}^j)}{x^j - w}
		\leq
		\mu_j \varepsilon_j
		\to
		0 .
	\end{equation*}
	Observing that $\mu_j \nabla_x \Lagr_{\mu_j}(\cdot,\widehat{y}^j)$ is locally Lipschitz continuous for all $\mu_j>0$ by \cref{ass:P}\eqref{ass:P:locLipschitzDiff},
	we have by \cite[Lemma 3.5]{demarchi2025mixed} and $x^j\to_J x^\star$ that $\{\Delta_j\}_{j\in J}$ remains bounded away from zero.
	Furthermore, using $\{\mu_j\}\searrow 0$ yields
	\begin{equation*}
		\mu_j \nabla_x \Lagr_{\mu_j}(x^j,\widehat{y}^j)
		\to_J
		\jac c(x^\star)^\top [c(x^\star) - \proj_{\CC}(c(x^\star))]
		=
		\nabla \feasmeas(x^\star)
	\end{equation*}
	by boundedness of $\{\widehat{y}^j\}$ and $\{\nabla f(x^j)\}_{j\in J}$, the latter due to $x^j\to_J x^\star$.	
	Overall, taking the limit $j\to_J\infty$, we have that
	\begin{align*}
		0
		={}&
		\lim_{j\to\infty}
		\max_{w\in\XX\cap\lpball(x^j,\Delta_j)}
		\innprod{\mu_j \nabla_x \Lagr_{\mu_j}(x^j,\widehat{y}^j)}{x^j - w} \\
		={}&
		\max_{w\in\XX\cap\lpball(x^\star,\Delta_\star)}
		\innprod{\nabla \feasmeas(x^\star)}{x^\star - w}
		={}
		\psimeas_{\feasmeas,\XX}(x^\star,\Delta_\star)
	\end{align*}
	for some $\Delta_\star>0$, proving the result.
\end{proof}

\subsection{Other Sequential Minimization Schemes}
\label{sec:OtherSchemes}

So far the focus has been on \cref{alg:AbstractALM},
but how do these developments affect other numerical approaches for \eqref{eq:P}?
Being part of the AL framework, the scheme analysed in \cite{grapiglia2020complexity} can be naturally extended to handle MINLPs.
Its peculiarity is that, starting with a feasible point, convergence to feasible accumulation points can be guaranteed, thanks to a reset mechanism.
Results similar to \cref{thm:globalConvergence,thm:minimalInfeasibility} can be readily obtained for this method too.
Indeed, analogous findings seem to extend far beyond the penalty scheme considered in \cref{sec:AugLagFramework},
possibly applying for a broad class of sequential minimization algorithms \cite{fiacco1968nonlinear}.
Although drawn in a different context, the arguments in \cite[Section~4]{demarchi2024interior} give a valid proof pattern for interior point (or barrier) methods, among others.

\newcommand{\ipval}{\mathcal{B}}
For illustrative purposes, let us consider the special case of \eqref{eq:P} with $\CC \coloneqq \spaceY_+$.
Introducing a barrier function $\func{b}{(0,\infty)}{\R}$ to approximate the indicator $\indicator_{\CC}$, e.g., the classical logarithmic barrier $b \colon t \mapsto - \log(t)$,
and a barrier parameter $\mu > 0$ to control this approximation,
one formulates a barrier subproblem---resembling \eqref{eq:ALsubproblem}---of the form
\begin{equation}
	\label{eq:ipsubproblem}
	\minimize~ \ipval_\mu(x) \quad
	\wrt~ x\in\XX 
	,\qquad
	\text{where}~
	\ipval_\mu(x)\coloneqq f(x) + \mu \sum_{i=1}^{m} b(c_i(x))
	.
\end{equation}
Then, a sequence of subproblems is solved, possibly inexactly and up to criticality, with decreasing barrier parameters.
This procedure is outlined in \cref{alg:AbstractIP}, where there is again no mention of $\Delta$.

Let us denote by $x^j$ an $\varepsilon_j$-critical point for the barrier subproblem \eqref{eq:ipsubproblem} with parameter $\mu_j > 0$.
Though with the drawback of requiring a strictly feasible point to start with (namely $x\in\XX$, $c(x)<0$),
at every iteration it must be that $x^j\in\XX$ and $c(x^j) < 0$,
that is, this barrier scheme maintains (strict) feasibility.
Moreover, echoing \cref{thm:globalConvergence}, it is easy to show that,
with $\mu_j, \varepsilon_j \to 0$,
accumulation points of $\{x^j\}$ are AKKT-critical for \eqref{eq:P};
see \cite[Thm~16]{demarchi2024interior}.
	Notice that the dual estimate rule at \cref{step:AbstractIP:y} of \cref{alg:AbstractIP} is justified by an identitiy analogous to \eqref{eq:gradLagrIdentity} for the augmented Lagrangian scheme, which now reads
\begin{equation}
	\nabla_x \ipval_\mu(x)
	={}
	\nabla f(x) + \jac c(x)^\top y
	={}
	\nabla_x \Lagr(x,y) .
\end{equation}
Finally, the update rule at \cref{step:AbstractIP:update_mu} forces the barrier parameter to vanish while remaining positive, so that the complementarity condition for KKT-criticality can be approximately satisfied; cf. \cite[Section 2]{waechter2006implementation} and \cite[Section 4]{demarchi2024interior}.
\begin{algorithm2e}[tbh]
	\DontPrintSemicolon
	Select $\mu_0, \varepsilon_0 > 0$, and
	$\kappa_\mu \in (0,1)$\;
	\For{$j = 0,1,2\ldots$}{
		Find an $\varepsilon_j$-critical point $x^j$ for $\ipval_{\mu_j}$ over $\XX$\label{step:AbstractIP:subproblem}\tcp*{subproblem}
		Set $y_i^j \gets \mu_j b^\prime(c_i(x^j))$, for all $i=1,\ldots,m$\label{step:AbstractIP:y}\;
		Set $\mu_{j+1} \gets \kappa_\mu \mu_j$\label{step:AbstractIP:update_mu}\;
		Select $\varepsilon_{j+1} \geq 0$ such that $\{\varepsilon_j\}\to 0$\label{step:AbstractIP:innertol}\;
	}
	\caption{Abstract interior point method for \eqref{eq:P}, with a barrier function $b$ suitable for $\CC$}
	\label{alg:AbstractIP}
\end{algorithm2e}

\section{Further Characterizations}
\label{sec:FurtherCharacterizations}

We now enrich the theoretical framework with results and interpretations well beyond those motivated by \cite{demarchi2025mixed} and \cref{alg:AbstractALM},
turning our attention to optimality conditions, Lagrangian duality, saddle point properties, and relationships with the classical proximal point algorithm.

For simplicity, 
we consider an optimization problem of the form \eqref{eq:P} with $\CC \coloneqq \KK$ a nonempty closed convex \emph{cone}.
Inspired by \cite[Section~8.4]{steck2018lagrange},
this assumption greatly simplifies the presentation thanks to the identity
\begin{equation}
	\label{eq:conj_polar_cone_identity}
	\indicator_{\KK}^\conj(y)
	=
	\sup_{z\in\KK} \innprod{z}{y}
	=
	\begin{cases}
		0 & \text{if}~y\in\KK^\polarcone, \\
		\infty & \text{otherwise}
	\end{cases}
	=
	\indicator_{\KK^\polarcone}(y)
	,
\end{equation}
which connects the indicator $\func{\indicator_{\KK}}{\spaceY}{\R\cup\{\infty\}}$ of a set $\KK\subseteq\spaceY$,
the \emphdef{conjugate function} $\func{h^\conj}{\spaceY}{\R\cup\{\infty\}}$ associated with a (proper and lower semicontinuous) function $\func{h}{\spaceY}{\R\cup\{\infty\}}$ \cite[Definition~13.1]{bauschke2017convex},
and the \emphdef{polar cone} $\KK^\polarcone\subseteq\spaceY$ of a subset $\KK$ of $\spaceY$ \cite[Definition~6.22]{bauschke2017convex}, respectively
\begin{equation*}
	h^\conj(v)
	\coloneqq
	\sup_{z\in\spaceY} \left\{ \innprod{z}{v} - h(z) \right\}
	\qquad\text{and}\qquad
	\KK^\polarcone
	\coloneqq
	\left\{ u\in\spaceY \,\middle\vert\, \sup_{v\in\KK} \innprod{v}{u} \leq 0 \right\}
	.
\end{equation*}

\subsection{Lagrangian Duality}
\label{sec:LagrangianDuality}

The necessary optimality conditions in \cref{def:KKTcritical} cannot be derived based on the Lagrangian function $\Lagr$ alone,
but additional insights on the problem are needed to setup the complementarity system encapsulated in the expression $y \in \normalcone_{\KK}(c(x))$.
Instead, a comprehensive first-order optimality analysis can be developed based on the \emph{generalized} Lagrangian function,
whose construction is briefly recalled
following \cite{rockafellar2023convergence,demarchi2024implicit,demarchi2024local}.
Introducing an auxiliary variable $s\in\spaceY$,
\eqref{eq:P} can be rewritten as
\begin{align}
	\tag{P$^\text{S}$}\label{eq:Pslack}
	\minimize~&f(x) \qquad\qquad
	\wrt~x\in\XX ,\, s\in \KK \\
	\stt~&c(x) - s = 0 , \nonumber
\end{align}
whose (classical) Lagrangian function, akin to \eqref{eq:classicalLagrangian}, reads
\begin{equation*}
	\Lagr^\text{S}(x,s,y)
	\coloneqq
	f(x) + \innprod{y}{c(x) - s}
	.
\end{equation*}
Marginalization of $\Lagr^\text{S}$ with respect to $s$ yields the generalized Lagrangian function $\func{\ell}{\XX\times\YY}{\R}$ associated to \eqref{eq:P}, given by
\begin{equation*}
	\ell(x,y)
	\coloneqq{}
	\inf_{s\in\KK} \Lagr^\text{S}(x,s,y)
	={}
	f(x) + \innprod{y}{c(x)} + \inf_s \left\{ \indicator_{\KK}(s) - \innprod{y}{s} \right\}
	={}
	\Lagr(x,y) - \indicator_{\KK}^\conj(y) .
\end{equation*}
Then,
observing the identity \eqref{eq:conj_polar_cone_identity},
the dual domain of $\ell$,
namely the set $\YY$ of valid multipliers, is given by
\begin{equation}
	\label{eq:dual_domain_polar_cone}
	\YY
	\coloneqq
	\spaceY\cap
	\dom \indicator_{\KK}^\conj
	=
	\dom \indicator_{\KK^\polarcone}
	=
	\KK^\polarcone
	,
\end{equation}
which corresponds to a nonempty closed convex cone in $\spaceY$.
Classical nonlinear programming is recovered by (neglecting integrality and) taking $\KK$ to be the standard constraint cone there:
$\KK\coloneqq\{0\}$ and $\KK\coloneqq\spaceY_-$
are associated respectively to $\YY\coloneqq\spaceY$ and $\YY\coloneqq\spaceY_+$.
Then, with this insight about the dual domain, a sound yet simple stratagem for providing a safeguarding set to \cref{alg:AbstractALM}
is to set $\Ybounded \coloneqq \YY \cap [-y_{\max},y_{\max}]^m$ for some large $y_{\max} > 0$ \cite[Section 3.1]{sopasakis2020open}.

In contrast with the (classical) Lagrangian $\Lagr$,
the emergence of dual information from the generalized Lagrangian $\ell$ allows not only to obtain dual estimates tailored to $\KK$,
but also to express primal-dual first-order optimality conditions without direct access to \eqref{eq:P}.
It is shown in \cite{rockafellar2023convergence}, \cite[Remark~3.5]{demarchi2024local}
that the generalized Lagrangian function $\ell$
is sufficient
to write necessary optimality conditions for \eqref{eq:P} when $\XX=\spaceX$ and $\KK$ is convex.
These read
\begin{equation}\label{eq:nocGeneralizedLagrangian}
	0 \in \partial_x \ell(x,y)
	\quad\text{and}\quad
	0 \in \partial_y (-\ell)(x,y),
\end{equation}
where the negative sign highlights the (generalized) saddle-point property of the primal-dual system.
But how does \eqref{eq:nocGeneralizedLagrangian} relate to \cref{def:KKTcritical}?
Owing to the identity $\nabla_x \Lagr = \nabla_x \ell$,
the first criticality condition $\psimeas_{\Lagr(\cdot,y),\XX}(x,\Delta)=0$ in \cref{def:KKTcritical} captures in fact
an extension of $0 = \nabla_x \ell(x,y)$ to accommodate the mixed-integer linear constraint set $\XX$.
Inspired by the descent-ascent motive behind \eqref{eq:nocGeneralizedLagrangian},
the main definition we will use below is the following,
with a character of primal-dual symmetry.

\begin{mydefbox}
	\begin{definition}\label[definition]{def:primal_dual_critical}
		A pair $(x,y) \in \XX\times \YY$ is called a \emphdef{local saddle point}
		of $\func{\Lagr}{\XX\times\YY}{\R}$ if
		\begin{equation*}
			\psimeas_{\Lagr(\cdot,y),\XX}(x,\Delta) = 0
			\quad\text{and}\quad
			\psimeas_{-\Lagr(x,\cdot),\YY}(y,\Delta) = 0
		\end{equation*}
		for some $\Delta>0$.
	\end{definition}
\end{mydefbox}
	Let us consider the set $\lpball(y,\Delta)$, which appears in the computation of $\psimeas_{-\Lagr(x,\cdot),\YY}(y,\Delta)$ according to \eqref{eq:psimeas}.
	Since $\YY$ is purely real-valued, the $\normlp{\cdot}$ norm there requires in fact no partial localization and therefore $\lpball(y,\Delta)$ is compact and convex.
	In this situation \cref{def:criticality} recovers classical criticality (or stationarity) notions for continuous optimization, for instance \cite[Definition 3.1]{byrd2005convergence}.
\begin{mybox}
	\begin{theorem}\label[theorem]{thm:local_saddle_Lagrangian}
		Consider \eqref{eq:P} and let $x\in\spaceX$, $y\in\spaceY$ be arbitrary but fixed.
		Then the following assertions are equivalent:
		\begin{enumerate}[label=(\roman{*})]
			\item $x$ is KKT-critical with multiplier $y$;
			\item $(x,y)$ is a local saddle point of $\Lagr$.
		\end{enumerate}
	\end{theorem}
\end{mybox}
\begin{proof}\label{proof:local_saddle_Lagrangian}
	Since both KKT-critical and local saddle points demand that $x\in\XX$ and $\psimeas_{\Lagr(\cdot,y),\XX}(x,\Delta) = 0$ holds for some $\Delta>0$,
	it remains to consider the second part of \cref{def:KKTcritical,def:primal_dual_critical}, namely the equivalence of
	$y\in\normalcone_{\KK}(c(x))$ and $\psimeas_{-\Lagr(x,\cdot),\YY}(y,\Delta) = 0$.
	We proceed by deriving a sequence of identities.
	Observing that
	\begin{equation*}
		0
		=
		\psimeas_{-\Lagr(x,\cdot),\YY}(y,\Delta)
		=
		\max_{w\in\YY\cap\lpball(y,\Delta)}
		\innprod{-\nabla_y\Lagr(x,y)}{y - w}
		\geq
		0
	\end{equation*}
	can be rewritten with a universival quantifier as
	\begin{equation*}
		\forall w\in\YY\cap\lpball(y,\Delta)
		\colon~
		\innprod{-\nabla_y\Lagr(x,y)}{y - w}
		=
		\innprod{y + \nabla_y\Lagr(x,y) - y}{w - y}
		\leq
		0 ,
	\end{equation*}
	the characterization \eqref{eq:projCharacterization} of projections onto convex sets yields
	\begin{equation*}
		y
		=
		\proj_{\YY \cap \lpball(y,\Delta)}\left( y + \nabla_y\Lagr(x,y) \right)
		.
	\end{equation*}
	Since all variables in $y$ are real-valued and the ball $\lpball(y,\Delta)$ is compact convex and centered at $y\in\YY$,
	the previous identity is equivalent to
	$
	y
	=
	\proj_{\YY}\left( y + \nabla_y\Lagr(x,y) \right)
	$
	for all $\Delta>0$.
	Using the property \eqref{eq:normalConeCharacterization} of normal cones
	and the partial derivative of $\Lagr$ in \eqref{eq:classicalLagrangian},
	we obtain
	$
	\nabla_y\Lagr(x,y)
	=
	c(x)
	\in
	\normalcone_{\YY}\left( y \right)
	$.
	Exploiting now the definition of $\YY$ \eqref{eq:dual_domain_polar_cone},
	the polar-conjugacy relation \eqref{eq:conj_polar_cone_identity} implies that
	$
	c(x)
	\in
	\partial\indicator_{\KK^\polarcone}(y)
	=
	\partial\indicator_{\KK}^\conj(y)
	$.
	Finally, owing to \cite[Proposition 11.3]{rockafellar2009variational},
	this is equivalent to
	$
	y\in\partial\indicator_{\KK}(c(x)) = \normalcone_{\KK}(c(x)) ,
	$
	which also implies the inclusion $c(x)\in\KK$,
	concluding the proof.
\end{proof}

\subsection{Saddle Points of the Augmented Lagrangian}

Inspired by the primal-dual characterization of KKT-critical points in \cref{sec:LagrangianDuality},
here we show that KKT-criticality for \eqref{eq:P} is also associated to a local saddle point property of the \emph{augmented} Lagrangian function.
This trait, recently re-investigated by Rockafellar \cite{rockafellar2023convergence} for a broad problem class,
allows to interpret the update rule at \cref{step:AbstractALM:y} as a dual gradient ascent step for the augmented Lagrangian,
thus making \cref{alg:AbstractALM} a primal descent, dual ascent method;
see also \cite[Section~8.1]{steck2018lagrange}.

We begin with some preliminary observations.
\begin{mybox}
	\begin{lemma}\label[lemma]{lem:auxResultPrimalFeas}
		Consider \eqref{eq:P} and let $x\in\XX$, $y\in\spaceY$, and $\Delta, \mu > 0$ be arbitrary but fixed.
		Then the following assertions are equivalent:
		\begin{enumerate}[label=(\roman{*})]
			\item\label{lem:auxResultPrimalFeas:normalCone}%
			$y\in\normalcone_{\KK}(c(x))$;
			\item\label{lem:auxResultPrimalFeas:gradLagr}%
			$\nabla_y \Lagr_\mu(x,y) = 0$;
			\item\label{lem:auxResultPrimalFeas:psiMeas}%
			$\psimeas_{-\Lagr_\mu(x,\cdot),\YY}(y,\Delta) = 0$.
		\end{enumerate}
		In particular, these conditions imply the inclusions $c(x)\in\KK$ and $y\in\YY$.
	\end{lemma}
\end{mybox}
\begin{proof}\label{proof:auxResultPrimalFeas}
	Owing to \eqref{eq:ALderivatives}, condition \ref{lem:auxResultPrimalFeas:gradLagr} can be rewritten as
	$c(x) = \proj_{\KK}(c(x) + \mu y)$
	and, since $\mu>0$, property \eqref{eq:normalConeCharacterization} implies the equivalence of \ref{lem:auxResultPrimalFeas:normalCone} and \ref{lem:auxResultPrimalFeas:gradLagr}.
	Now, patterning the proof of \cref{thm:local_saddle_Lagrangian}, we obtain that \ref{lem:auxResultPrimalFeas:psiMeas} is equivalent to
	$
	\nabla_y \Lagr_\mu(x,y)
	\in
	\normalcone_{\YY}\left( y \right)
	$.
	Then,
	the implication \ref{lem:auxResultPrimalFeas:gradLagr}$\implies$\ref{lem:auxResultPrimalFeas:psiMeas} is clear,
	and it remains to focus on the converse one.
	
	Let us consider now the maximization of $\Lagr_\mu(x,\cdot)$ over $\spaceY$, that is, dropping the restriction to $\YY$---as well as the trust region in \eqref{eq:psimeas}.
	Then, any (unconstrained) solution $\widetilde{y}\in\spaceY$ necessarily satisfies $\nabla_y \Lagr_\mu(x,\widetilde{y}) = 0$,
	which is equivalent to
	$\widetilde{y} \in \normalcone_{\KK}(c(x))$ by combining \eqref{eq:ALderivatives}--\eqref{eq:ALauxiliaries} and \eqref{eq:normalConeCharacterization}.
	Furthermore, owing to convexity of $\KK$ and \cite[Proposition~11.3]{rockafellar2009variational},
	this inclusion coincides with
	$c(x) \in \normalcone_{\KK^\polarcone}(\widetilde{y})$,
	meaning in particular that $\widetilde{y}\in\KK^\polarcone=\YY$ by \eqref{eq:dual_domain_polar_cone}.
	Thus, since the unconstrained optimum $\widetilde{y}$ satisfies in fact the restriction to $\YY$, 
	it is optimal for the constrained problem too.
	Indeed, by convexity of $\YY$, $\widetilde{y}$ remains optimal also considering a trust region $\lpball(\widetilde{y},\Delta)$, for any $\Delta>0$,
	thus showing that \ref{lem:auxResultPrimalFeas:psiMeas}$\implies$\ref{lem:auxResultPrimalFeas:gradLagr}.
	
	Finally, the inclusions follow respectively from the normal cone $\normalcone_{\KK}(c(x))$ being nonempty in \ref{lem:auxResultPrimalFeas:normalCone} and from the restriction $y\in\YY$ in \eqref{eq:psimeas} for \ref{lem:auxResultPrimalFeas:psiMeas}.
\end{proof}

The following is the main result of this section.
\begin{mybox}
	\begin{theorem}\label[theorem]{thm:saddlePointAugmentedLagrangian}
		Consider \eqref{eq:P} and let $x\in\spaceX$, $y \in \spaceY$ be arbitrary but fixed.
		Then the following assertions are equivalent:
		\begin{enumerate}[label=(\roman{*})]
			\item $x$ is KKT-critical with multiplier $y$;
			\item $(x,y)$ is a local saddle point of $\Lagr_\mu$ for some $\mu>0$;
			\item $(x,y)$ is a local saddle point of $\Lagr_\mu$ for all $\mu>0$.
		\end{enumerate}
	\end{theorem}
\end{mybox}
\begin{proof}\label{proof:saddlePointAugmentedLagrangian}
	We prove the equivalence via a loop of implications.
	Note that $(iii) \implies (ii)$ is straightforward.
	
	For the implication $(ii) \implies (i)$,
	let $(x,y)$ be a local saddle point of $\Lagr_\mu$ for some $\mu>0$.
	Then \cref{lem:auxResultPrimalFeas} implies that $c(x) \in \KK$ and $y \in \normalcone_{\KK}(c(x))$. 
	Therefore, by combining with \eqref{eq:ALderivatives}--\eqref{eq:ALauxiliaries} and properties \eqref{eq:projCharacterization}--\eqref{eq:normalConeCharacterization}, we obtain the identity
	\begin{equation}
		\label{eq:dxAL_equal_dxL}
		\nabla_x \Lagr_\mu(x,y)
		={}
		\nabla f(x) + \jac c(x)^\top y
		={}
		\nabla_x \Lagr(x,y).
	\end{equation}
	Therefore, since $\psimeas_{\Lagr_\mu(\cdot,y),\XX}(x,\Delta) = 0$ holds for some $\Delta>0$, it must be also $\psimeas_{\Lagr(\cdot,y),\XX}(x,\Delta) = 0$.
	Thus, $x$ is KKT-critical for \eqref{eq:P} with multiplier $y$.
	
	For the remaining implication $(i)\implies(iii)$,
	let $\mu>0$ be arbitrary but fixed and $x$ a KKT-critical point with multiplier $y$.
	Then, $c(x)\in\KK$ and $y \in \normalcone_{\KK}(c(x))$ hold owing to KKT-criticality.
	Hence, on the one hand, \cref{lem:auxResultPrimalFeas} implies that the second equality in \cref{def:primal_dual_critical} is satisfied.
	On the other hand, this furnishes again \eqref{eq:dxAL_equal_dxL}, and thus KKT-criticality of $(x,y)$ yields $\psimeas_{\Lagr_\mu(\cdot,y),\XX}(x,\Delta) = \psimeas_{\Lagr(\cdot,y),\XX}(x,\Delta) = 0$.
	With $\mu>0$ being arbitrary, this shows that $(x,y)$ is a local saddle point of $\Lagr_\mu$ for all $\mu >0$.
\end{proof}

\subsection{Relationship with Proximal Point Methods}

Connections of augmented Lagrangian methods with duality and the proximal point algorithm (PPA) have been discussed in Hilbert spaces \cite[Section~8.4]{steck2018lagrange} and explored in the broad setting of generalized nonlinear programming \cite{rockafellar2023convergence,demarchi2024local}.
We turn now to examining these properties in the context of MINLP.
Considering \eqref{eq:P}, the associated Lagrangian function \eqref{eq:classicalLagrangian}, and the dual domain $\YY$ \eqref{eq:dual_domain_polar_cone},
we define for all $y\in\YY$
\begin{equation*}
	\dualLagr(y)
	\coloneqq
	\inf_{x\in\XX} \Lagr(x,y)
	=
	\inf_{x\in\XX} \left\{ f(x) + \innprod{y}{c(x)} \right\}
\end{equation*}
so that the natural ``dual'' problem of \eqref{eq:P} is given by
\begin{equation*}
	\maximize~ \dualLagr(y) \quad
	\wrt~ y\in\YY
	.
\end{equation*}
Note that $\dualLagr$ is a concave function since it is an infimum of affine functions.
Then, by convexity of $\YY$, the above is a concave maximization problem,
equivalent to a convex minimization problem.
Given a starting point $y^0$,
the PPA consists in applying the recursion
\begin{equation*}
	y^{j+1} \coloneqq \prox_{-\nu_j \dualLagr}(y^j)
\end{equation*}
with parameter $\nu_j > 0$,
where the central ingredient is the \emphdef{proximal mapping} associated to the problem, given by
\begin{equation*}
	\prox_{- \nu \dualLagr}(w)
	\coloneqq
	\argmin_{y\in\YY} \left\{ - \dualLagr(y) + \frac{1}{2\nu} \| y - w \|^2 \right\}
\end{equation*}
for any $\nu > 0$ \cite[Chapter~24]{bauschke2017convex}, \cite[Section~2]{rockafellar2023convergence}.
Note that the function occurring inside the $\argmin$ is strongly convex,
hence it admits a unique minimizer, and thus the proximal mapping is well-defined and single-valued.
We will demonstrate that this iterative procedure is (still) strongly related to the AL method, whose basic iteration with parameter $\mu_j > 0$ reads
\begin{equation*}
	x^{j+1} \in{} \argmin_{x\in\XX} \Lagr_{\mu_j}(x,y^j) ,\quad
	z^{j+1} \coloneqq{} \proj_{\KK}(c(x^{j+1}) + \mu_j y^j) ,\quad
	y^{j+1} \coloneqq{} y^j + \frac{c(x^{j+1}) - z^{j+1}}{\mu_j},
\end{equation*}
where $z^{j+1}\in\KK$ and $y^{j+1}\in \normalcone_{\KK}( z^{j+1} )$ hold by construction; see \cref{lem:iterates}.

The main result in this section is the following \cref{thm:proxPointMethod0}, which shows that, up to criticality, a basic AL method for \eqref{eq:P} is equivalent to applying PPA to the dual problem.

\begin{mybox}
	\begin{theorem}\label[theorem]{thm:proxPointMethod0}
		Consider \eqref{eq:P} and let $w \in \spaceY$, $\mu > 0$ be arbitrary but fixed.
		Let $\bar{x}$ be a critical point for $\Lagr_\mu(\cdot,w)$ over $\XX$.
		Define
		$\bar{s} \coloneqq \proj_{\KK}(c(\bar{x}) + \mu w)$ and
		$\bar{y} \coloneqq w + [c(\bar{x}) - \bar{s}]/\mu$.
		Then $\bar{y} = \prox_{-\mu \dualLagr}(w)\in\YY$ and $\bar{x}\in\XX$ is a critical point for the infimum defining $\dualLagr(\bar{y})$, namely for $\Lagr(\cdot,\bar{y})$ over $\XX$.
	\end{theorem}
\end{mybox}
\begin{proof}\label{proof:proxPointMethod0}
	We prove the claim by showing that $(\bar{x},\bar{y})$ is a local saddle point of the function
	\begin{equation*}
		\func{h}{\XX \times \YY}{\R}
		,\qquad
		h(x,y) \coloneqq \Lagr(x,y) - \frac{\mu}{2} \| y- w \|^2
	\end{equation*}
	which brings together the dual function $\dualLagr$ with the quadratic proximal term.
	To verify this saddle property, note that the definition of $\bar{x}$ and $\bar{y}$ implies by \eqref{eq:ALderivatives}--\eqref{eq:ALauxiliaries} that
	\begin{equation*}
		\nabla_x \Lagr_\mu(\bar{x},w)
		=
		\nabla f(\bar{x}) + \jac c(\bar{x})^\top \bar{y}
		=
		\nabla_x \Lagr(\bar{x},\bar{y})
		=
		\nabla_x h(\bar{x},\bar{y}).
	\end{equation*}
	Then, by \cref{def:criticality}, there exists some $\Delta>0$ such that
	\begin{equation*}
		0
		=
		\psimeas_{\Lagr_\mu(\cdot,w),\XX}(\bar{x},\Delta)
		=
		\psimeas_{\Lagr(\cdot,\bar{y}),\XX}(\bar{x},\Delta)
		=
		\psimeas_{h(\cdot,\bar{y}),\XX}(\bar{x},\Delta)
		,
	\end{equation*}
	hence $\bar{x}$ is a critical point for
	$h(\cdot,\bar{y})$ over $\XX$.
	On the other hand, $h(\bar{x},\cdot)$ is a strictly concave quadratic function of the form
	\begin{equation*}
		h(\bar{x},\cdot)
		\colon
		y
		\mapsto
		- \frac{\mu}{2} \left\| y - w + \frac{c(\bar{x})}{\mu} \right\|^2 + c_h ,
	\end{equation*}
	where $c_h\in\R$ is a constant independent of $y$.
	Therefore, the unique maximizer $\widetilde{y}$ of $h(\bar{x},\cdot)$ over the convex set $\YY$ is
	determined by the necessary optimality condition $\nabla_y h(\bar{x},\widetilde{y}) \in \normalcone_{\YY}(\widetilde{y})$.
	Using the definition of $h$, \eqref{eq:ALderivatives}--\eqref{eq:ALauxiliaries}, \eqref{eq:dual_domain_polar_cone}, and the identity \eqref{eq:conj_polar_cone_identity},
	this can be rewritten as
	\[
	c(\bar{x}) + \mu (w - \widetilde{y})
	\in
	\normalcone_{\KK^\polarcone}(\widetilde{y})
	=
	\partial\indicator_{\KK}^\conj(\widetilde{y})
	.
	\]
	Then, by convexity of $\KK$ and \cite[Proposition~11.3]{rockafellar2009variational},
	this is equivalent to
	$
	\widetilde{y} \in
	\normalcone_{\KK}(c(\bar{x}) + \mu (w - \widetilde{y}))
	$.
	Finally, the definition of $\bar{s}$ and characterization \eqref{eq:normalConeCharacterization} yield the identity
	\begin{equation*}
		\bar{s} \coloneqq \proj_{\KK}(c(\bar{x}) + \mu w) = c(\bar{x}) + \mu (w - \widetilde{y})
		,
	\end{equation*}
	showing that the unique maximizer $\widetilde{y}$ coincides in fact with $\bar{y}$,
	concluding the proof.
\end{proof}

\section{Concluding Remarks}

The developments and results in this paper offer solid theoretical foundations for employing continuous optimization techniques to address mixed-integer nonlinear programming,
at least as principled heuristics.
Although presented in details for an augmented Lagrangian scheme, a similar analysis readily applies to other sequential minimization techniques, such as barrier and mixed schemes.
Preliminary numerical tests on the optimal control of hybrid dynamics demonstrated the viability of the proposed approach, but
only a more comprehensive computational validation and comparison will attest its practical performance, limitations, and range of applications.
We foresee the need for combining solvers to deliver, exploiting warm-starts, good quality solutions with low computational effort.

It remains an open question how to relax the requirements on the problem data, particularly \cref{ass:P}\eqref{ass:P:integerBounded}, which however concerns the subsolver only.
When localizing both real- and integer-valued variables, enough freedom should be left for the latter, but not necessarily for the former.
In particular, one should prevent that some integers become effectively fixed,
leading to weaker optimality conditions.

\ifpreprint
	\phantomsection
	\addcontentsline{toc}{section}{References}%
	{\small
		\bibliographystyle{habbrv}
		\bibliography{biblio}
	}
\fi

\appendix
\section{Motivating Numerical Example}\label{sec:example_mila_cia}

This appendix illustrates with a numerical example some of the benefits that come with the mixed-integer linearization scheme of \cite{demarchi2025mixed}.
The problem \eqref{eq:bocp} below is in fact a MIQP that could be solved by specialized solvers in a matter of milliseconds.
However, the purpose of this section is to show that, even for the toy problem \eqref{eq:bocp}, the CIA approach returns a suboptimal solution that can be significantly improved by the integrality-preserving MILA of \cite{demarchi2025mixed}.
Nonetheless, CIA is a scalable approach that often generates good initial approximations for further refinement with MILA.

\medskip

Let us consider the optimal control of a discrete-time linear dynamics with binary-valued control, with one state, one control input, and quadratic tracking cost.
Combinatorial constraints are incorporated in the form of a maximum number of switches for the control input.
The problem formulation reads
\renewcommand{\xi}{s}
\begin{subequations}
	\label{eq:bocp}
	\begin{align}
		\minimize~&h \sum_{k=0}^N (\xi_k-1)^2 &
		\wrt~&\{\xi_k\}_{k=0}^N ,\, \{b_k\}_{k=0}^{N-1} \label{eq:bocp_cost}\\
		\stt~&\xi_{k+1} = \xi_k + h (b_k - \tfrac{1}{2} ) ,\quad b_k \in \{0,1\} & \text{for}~&k=0,\ldots,N-1, \label{eq:bocp_dynamics}\\
		&\xi_0 = 0 = \xi_N , \label{eq:bocp_bc}\\
		&\sum_{k=0}^{N-2} |b_{k+1} - b_k| \leq \sigma_{\max} , \label{eq:bocp_switch}
	\end{align}
\end{subequations}
where $h\coloneqq T/N$ is the time step, with $T\coloneqq 10$ and $N\coloneqq 100$, 
$\xi_k$ and $b_k$ denote the discrete-time state and control, respectively, at time $t_k = k h$, $k\in\N$.
The objective function in \eqref{eq:bocp_cost} promotes state values near one, while initial and terminal conditions in \eqref{eq:bocp_bc} require the state to be zero there.
As the control is binary-valued, the dynamics in \eqref{eq:bocp_dynamics} prevent the state from remaining constant.
The summation term in the inequality constraint in \eqref{eq:bocp_switch} counts the number of switches, namely how many times the control input changes value in $\{0,1\}$.
The maximum number of switches allowed is $\sigma_{\max}\coloneqq10$.
It should be noted that, since the absolute value can be recast into linear inequalities at the price of some auxiliary variables, all constraints in \eqref{eq:bocp} can be written in mixed-integer linear form.

\begin{figure}[tbh]
\centering%
\includegraphics{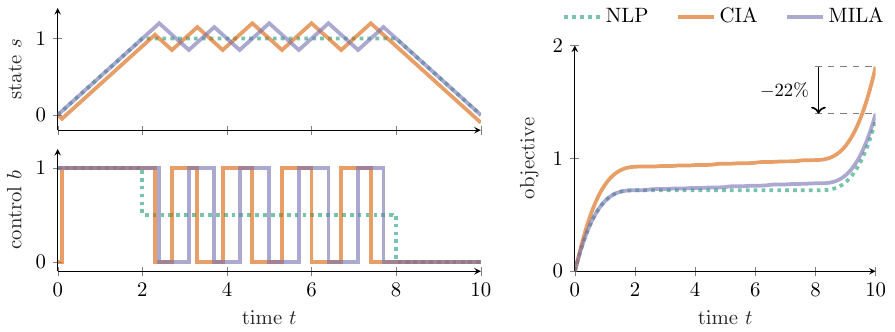}
\caption{%
	Results for the binary optimal control problem \eqref{eq:bocp} with $N=100$ discretization intervals: solutions obtained with relaxed integrality
	(NLP), combinatorial integral approximation (CIA), and (warm-started) mixed-integer linearization algorithm (MILA).%
}%
\label{fig:bocp}%
\end{figure}

The first step of \cite{sager2011combinatorial}'s decomposition method is to relax the integrality constraint in \eqref{eq:bocp}, replacing $\{0,1\}$ with $[0,1]$, and solve the corresponding NLP (convex in this case).
The relaxed solution obtained with \texttt{Ipopt}%
\footnote{%
	Version 3.14.16, with the option \texttt{tol} set to $10^{-8}$ and \texttt{honor\_original\_bounds} to \texttt{yes}.
}
\cite{waechter2006implementation} is depicted in \cref{fig:bocp} (labelled ``NLP'').
After an initial phase the control settles around the optimal value $\nicefrac{1}{2}$, for which the state can track exactly one and the overall cost $J_{\text{NLP}}\approx 1.435$ is a lower bound for binary control strategies.
Although solved without switching constraint, the relaxed control input switches only twice, and therefore it is feasible for \eqref{eq:bocp}.

The second step is the so called \emph{combinatorial integral approximation} (CIA):
starting from the relaxed control input, a binary-valued sequence is obtained from the software package \texttt{pycombina}%
\footnote{%
	Version 0.3.4, using the tailored \texttt{CombinaBnB} solver with the option \texttt{max\_iter} set to $10^9$.
}
\cite{buerger2020pycombina} with an explicit specification of the switching constraint.
The ``CIA'' solution is also depicted in \cref{fig:bocp}, exhibiting exactly $\sigma_{\max}$ switches and an increased cost $J_{\text{CIA}}\approx 1.934$ due to degraded tracking performance.
Moreover, the CIA solution does not satisfy the terminal condition.

Finally, we adopt the mixed-integer linear algorithm (MILA) of \cite{demarchi2025mixed}, which takes into account both the system dynamics and the combinatorial constraints.
Using the CIA solution as starting point, Algorithm~3.1 of \cite{demarchi2025mixed}%
\footnote{%
	Version 0.1.5, with monotone decrease and tolerance \texttt{neg\_tol} for negative criticality values set to $10^{-14}$.
}
generates feasible iterates with improved cost.
The solver returns after 5 iterations with cost $J_{\text{MILA}}\approx 1.5035$, with a dramatic $-22\%$ cost reduction relative to $J_{\text{CIA}}$, which brings the MILA solution to be only 4.8\% above the (unattainable) $J_{\text{NLP}}$ lower bound.

This simple example demonstrates that MILA can improve upon the solutions delivered by the state-of-the-art decomposition method \cite{sager2011combinatorial}.
However, it cannot be stressed enough that good quality local solutions can be achieved in reasonable time only by combining (and warm-starting) these techniques.

\section{Numerical Experience with Nonlinear Constraints}
\label{sec:NumericalResults}

This appendix is dedicated to testing the proposed framework and to assessing its numerical performance.
The example problem detailed below concerns the optimal point-to-point control of a point-mass with hybrid dynamics.
The model captures the (nonlinear) longitudinal dynamics of a car with aerodynamic drag and downforce, while a turbo charger mechanism gives rise to mixed-integer linear constraints.
A \texttt{python} implementation of \cref{alg:AbstractALM,alg:AbstractIP}, denoted respectively ``AL'' and ``IP'', is invoked with different model parameters, discretizations and initial guesses.
In particular, due to the tradeoffs in affordable optimization methods, we observe different regimes when using an all-zero initialization and one obtained through a simple heuristic.
For comparison, a \texttt{C++} implementation of dynamic programming (``DP'') is considered as baseline, since it does not require an initial guess and its results are globally optimal (up to the discretization).

\paragraph*{Implementation details}
AL and IP rely on a \texttt{python} implementation of MILA \cite[Alg. 3.1]{demarchi2025mixed} as subsolver, which in turn calls \texttt{Gurobi} (version 11.0.0) as MILP solver.
The IP implementation handles nonlinear inequality constraints via the logarithmic barrier function $b \coloneqq -\log$; equality constraints are treated as in AL.
The algorithmic parameters in \cref{alg:AbstractALM,alg:AbstractIP} are set as follows:
$\varepsilon_0 = \mu_0 = 0.1$, $\kappa_\mu = 0.5$, $\theta_\mu=0.9$, $\eta_j=\epsilon$ and $\varepsilon_{j+1} = \max\{ \epsilon,\kappa_\varepsilon\varepsilon_j \}$ for all $j\in\N$ with $\kappa_\varepsilon = 0.5$ and termination tolerance $\epsilon$.
The dual safeguarding set $\Ybounded$ is a hyperbox, defined by $[-y_{\max},y_{\max}]$ for each equality $c_i(x)=0$ and $[0,y_{\max}]$ for each inequality $c_i(x)\leq 0$, with $y_{\max}\coloneqq 10^{20}$.
AL and IP terminate and return $(x^j,y^j)$ when $\epsilon$-KKT-criticality is detected.
We set the tolerance $\epsilon \coloneqq 10^{-6}$.

\subsection*{Illustrative Problem}

The problem under consideration extends the optimal control example in \cite[Section 4.1]{nikitina2025hybrid}, including nonlinear dynamics, integrality requirements, mixed state-control constraints and path inequality specifications.
The double-integrator point-mass model of a car is equipped with a hysteretic turbo accelerator.
More specifically, the car's state is described by its position $s(t)$, velocity $v(t)$ and turbo state $w(t)\in\{0,1\}$, which are governed by
\[
\dot{s}(t) = v(t) ,\qquad
\dot{v}(t) = \tau(w(t),a(t)) - b(t) - c_d v(t)^2 ,
\]
and by the hysteresis curve: the turbo mode becomes active ($w=1$) when the velocity exceeds $v^+ \coloneqq 10$ and it becomes inactive ($w=0$) when the velocity falls below $v^- \coloneqq 5$.
The velocity is limited by $|v(t)| \leq v_{\max} \coloneqq 25$.
The car is controlled with the input to the acceleration and brake pedals, respectively $a(t)\in[0,a_{\max}]$ and $b(t)\in [0,b_{\max}]$, with $a_{\max}\coloneqq 5$ and $b_{\max}\coloneqq 10$.
The traction $\tau$ has two modes of operation depending on the turbo state, defined as
$\tau(w,a)=a$ if $w=0$, and $\tau(w,a)=3a$ if $w=1$.
Parameter $c_d \coloneqq 10^{-3}$ denotes the drag coefficient.

Given the final time $T\coloneqq 10$, the task is to bring the car from the initial state $(s(0),v(0))=(0,0)$ to the final state $(s(T),v(T))=(150,0)$ with minimum effort, as encoded by the objective
\[
\min \int_{0}^{T} [a(t)^2 + \alpha_b b(t)^3] \mathrm{d}t ,
\]
where $\alpha_b\coloneqq 10^{-2}$.
Finally, we model a limitation of grip in the form of bilateral path constraints, requiring that the tangential force does not exceed a certain fraction of the normal force between car and road surface, namely that
\[
| \tau(w,a) - b | \leq c_z + c_g v^2
\]
holds, where parameters $c_z > 0$ and $c_g\coloneqq 10^{-3}$ identify the grip quality (low values correspond to low grip).
This additional (nonlinear inequality) constraint in the model gives us the opportunity to showcase and compare the AL and IP strategies.

\paragraph*{CIA and DP}
The hybrid turbo dynamics is difficult to formulate in partial outer convexification form, if possible at all, hindering the application of CIA for systems with state-dependent jumps \cite{sager2011combinatorial}.
Moreover, mixed state-control constraints are not included in the binary reconstruction step of the original CIA; see \cite{zeile2021mixed,buerger2020pycombina,buerger2023gauss} for some recent developments.
Conversely, the application of DP on the (discretized) hybrid dynamics is straightforward, but path constraints and final state conditions are not easily incorporated and must be treated with penalty terms.
The violation of final conditions is penalized as a Mayer term, namely adding to the objective the cost
\[
\lambda_{\rm DP}[(s(T)-150)^2 + v^2(T)]
\]
with $\lambda_{\rm DP}\coloneqq100$.
Analogously, the grip constraint is incorporated as a Lagrange cost of the form
\[
\lambda_{\rm DP} \int_{0}^T \max\{0, |\tau(a(t),w(t)) - b(t)| - c_z - c_g v^2(t) \}^2 \mathrm{d}t
.
\]
Finally, in addition to the time discretization, dynamic programming requires state and control grids:
the position is discretized with $100$ intervals over the range $[0,150]$, the velocity with $50$ over $[0,25]$, the turbo state is binary, the acceleration and brake pedals with $20$ intervals over $[0,5]$ and $[0,10]$ respectively.
The selected discretization and penalty parameter $\lambda_{\rm DP}$ strike a balance between errors in final position and velocity (less than 1) and manageable runtimes.

\paragraph*{Time discretization}
The optimal control problem is cast in the form of \eqref{eq:P} by introducing a time grid with $N$ intervals over $[0,T]$.
Adopting the explicit Euler scheme, the dynamics of $s$ and $v$ become a set of $2N$ equality constraints.
Then, the finite-dimensional model has $3(2N+1)$ variables: $2(N+1)$ for the real-valued states $s$ and $v$, $N+1$ binary-valued for $w$, $2N$ for the controls $a$ and $b$, $N$ for the auxiliary $\tau$.
The grip constraint leads to $2N$ nonlinear inequalities, which are treated with either a shifted penalty (AL) or a barrier (IP) approach.
The logical implications describing the hysteresis characteristic are specified by $8N$ mixed-integer linear constraints, as detailed in \cite[Section 4.1]{demarchi2025mixed}.
Numerical results below are presented up to $N=100$, which corresponds to hundreds of (real and binary) variables and (linear and nonlinear) constraints.

\paragraph*{Initial guess}
The AL and IP solvers will be invoked with two kinds of initial guesses, with the goal of inspecting their behaviour in different circumstances.
An all-zero initialization simulates a cold-start for the solver, as it is relatively far from an optimal solution.
In contrast, an improved initialization provides a warm-start for the solver.
This is obtained by integrating the discrete-time dynamics with heuristic control inputs:
first, 90\% of the maximum acceleration pedal is applied until 90\% of the speed limit is reached, then the acceleration is graded to maintain this constant speed before applying 90\% of the maximum brake pedal to reach zero velocity at the final time $T$.

\subsection*{Results and Comparisons}

\begin{figure}[tbh]
	\centering
	\includegraphics{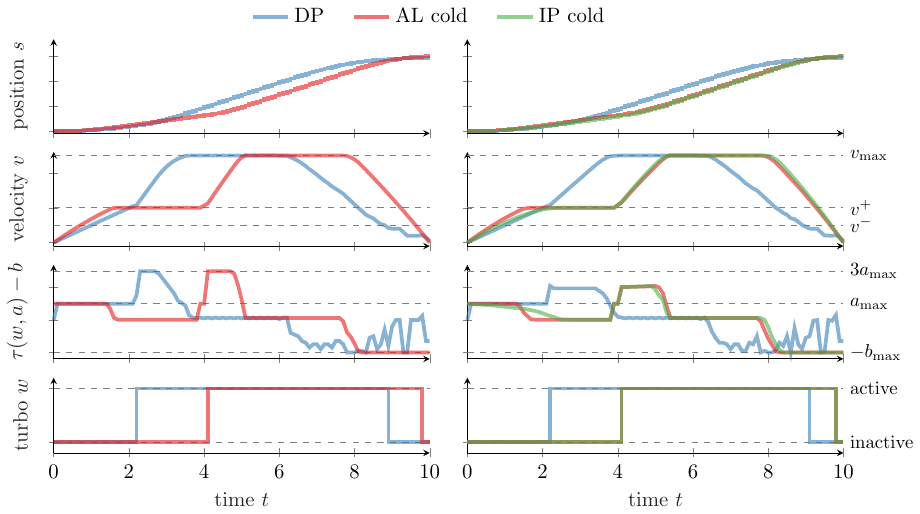}
	\caption{%
		Results of the turbo car problem discretized with $N=100$, for $c_z=\infty$ (left) and $c_z=10$ (right).
		Comparison of cold-started AL and IP, starting from an all-zero initial guess, against DP.%
	}%
	\label{fig:turbocar_grip_cold}
\end{figure}

\begin{figure}[tbh]
	\centering
	\includegraphics{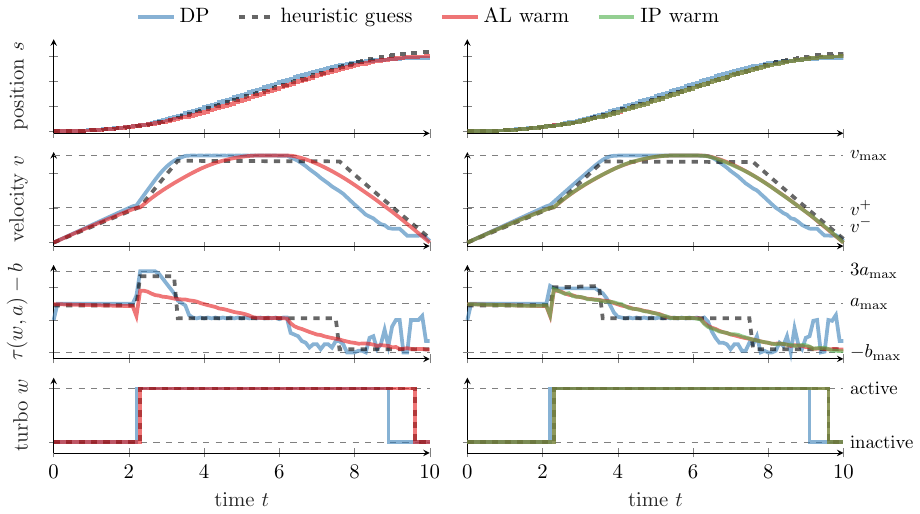}
	\caption{%
		Results of the turbo car problem discretized with $N=100$, for $c_z=\infty$ (left) and $c_z=10$ (right).
		Comparison of warm-started AL and IP, starting from a heuristic initial guess, against DP.%
	}%
	\label{fig:turbocar_grip_warm}
\end{figure}

The solutions returned by AL, IP and DP are depicted in \cref{fig:turbocar_grip_cold,fig:turbocar_grip_warm}, respectively with cold- and warm-starting.
Since the grip constraint does not apply when $c_z=\infty$, the IP strategy appears only for the case $c_z=10$.
The DP solution recovers the optimal pattern found in \cite{demarchi2025mixed,nikitina2025hybrid}, but the controls exhibit additional oscillations in the final section; these artifacts are likely due to the state and control discretization.
The cold-started AL and IP return the same feasible but possibly suboptimal trajectory: compared to the DP solution, the turbo activation is delayed and the final phase requires maximum braking, which is uncommon for a minimum-effort control task.
When warm-started with the heuristic initial guess, AL and IP generate feasible trajectories with a turbo activation pattern closer to the DP solution, as shown in \cref{fig:turbocar_grip_warm}, and with much smoother control inputs.

\begin{figure}[tbh]
	\centering
	\includegraphics{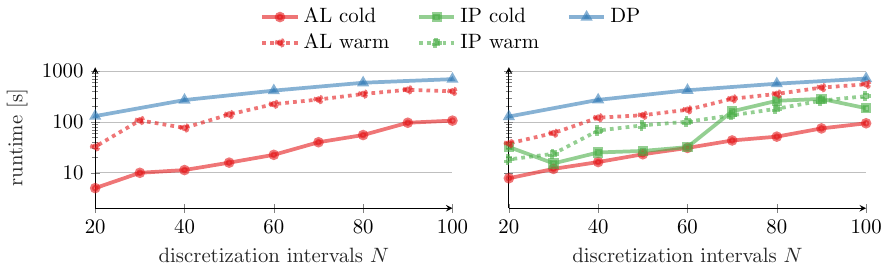}
	\caption{%
		Runtimes for the turbo car problem discretized with different number $N$ of intervals.
		Results obtained for $c_z=\infty$ (left) and $c_z=10$ (right), with all-zero (cold) and heuristic (warm) initial guesses.
		DP requires no starting point.%
	}%
	\label{fig:turbocar_runtime}
\end{figure}

The proposed affordable solvers are compared also based on their runtimes, which are summarized in \cref{fig:turbocar_runtime} for $N\in \{20,40,\ldots,100\}$.
The computational effort grows linearly with $N$ for DP and faster for AL and IP.
Nevertheless, DP takes the longest runtime on each instance (and requires a large working memory), despite the coarse (state and control) discretization and the parallelization of execution on 12 cores.
Conversely, the performance of AL and IP can strongly depend on the initial guess provided, as highlighted by the consistent and considerable difference between cold- and warm-started executions.
This feature is typical of affordable methods, as the requirement of global optimality is relaxed, seeking a tradeoff between solution quality and computational effort.

The results in \cref{fig:turbocar_runtime} together with \cref{fig:turbocar_grip_cold,fig:turbocar_grip_warm} can be interpreted as follows:
when cold-started, the iterates quickly accumulate at a local minimizer with a simple (almost piecewise constant) control sequence;
when warm-started, the iterates approach a more complicated, higher-quality control sequence which requires refinement, and so more iterations.
This sensitivity does not affect the global DP approach, which explores the whole state-control space and uses no initial guess.
In contrast, since DP relies on state and control grids while AL and IP do not, the solution obtained from the latter solvers can be much more accurate, as demonstrated by the low termination tolerance $\epsilon=10^{-6}$ compared to the coarse discretization for DP.
Moreover, even though AL and IP adopt a first-order inner solver, namely MILA from \cite{demarchi2025mixed}, which exhibits a slow tail convergence, their runtimes are still better than DP's, and with a reduced memory footprint.

Overall, this preliminary numerical investigation showcases not only the potential of the proposed mixed-integer Lagrangian framework in applications, but also the modeling flexibility it offers.

\end{document}